\documentclass[11pt,handout]{article}%
\usepackage{geometry}                % See geometry.pdf to learn the layout options. There are lots.
\usepackage{tikz} 
\usetikzlibrary{positioning}
\usepackage{graphicx}
\usepackage{amssymb}
\usepackage{epstopdf}
\usepackage{subfigure}  % use for side-by-side figures
\usepackage[sans]{dsfont}
\usepackage[english]{babel}
\usepackage{latexsym}
\usepackage{mathrsfs}
\usepackage{graphicx}
\usepackage{color}
\usepackage{float}
\usepackage{mathtools}
\usepackage{amsmath}
\usepackage{amsfonts}
\numberwithin{equation}{section}
\usepackage{enumerate}
\usepackage{amsthm}
\usepackage[title]{appendix}

\usepackage[bookmarks=true,colorlinks=true,linkcolor={blue},urlcolor={blue}, citecolor={blue},pdfstartview={XYZ null null 1.22}]{hyperref}%
\def\be#1\ee{\begin{equation}#1\end{equation}}

\newtheorem{proposition}{Proposition}
\theoremstyle{definition}

\newtheorem{remark}{Remark}

\newcommand{\xx}{{\bf x}}
\newcommand{\rr}{{\bf r}}
\newcommand{\vvel}{{\bf v}}
\newcommand{\BB}{{\bf B}}
\newcommand{\EE}{{\bf E}}
\def\RR{\mathbb R}

\title{Control strategies for magnetized plasma: a polar coordinates framework}

\author{ 
	F. Ferrarese \thanks{Department of Mathematics and Computer Science \& Center for Modeling, Computing and Statistics (CMCS), University of Ferrara, via Machiavelli 30, 44121 Ferrara, ITALY (federica.ferrarese@unife.it).}  
	}
 
\date{}
\begin{document}
	\maketitle
	
	\section*{Abstract}
	In this work, we provide an overview of various control strategies aimed at steering plasma toward desired configurations using an external magnetic field. From a modeling perspective, we focus on the Vlasov equation in a two-dimensional bounded domain, accounting for both a self-induced electric field and a strong external magnetic field. The results are presented in a polar coordinate framework, which is particularly well-suited for simulating toroidal devices such as Tokamaks and Stellarators. A key feature of the proposed control strategies is their feedback mechanism, which is based on an instantaneous prediction of the discretized system. Finally, different numerical experiments in the two-dimensional polar coordinate setting demonstrate the effectiveness of the approaches.

\section{Introduction} 
In recent years, considerable attention has been devoted to developing advanced numerical techniques for simulating plasma physics phenomena \cite{cheng1976integration,filbet2003numerical,filbet2018numerical,ghizzo1993eulerian, dimarco2015numerical}. Of particular interest are magnetized plasmas, which are central to nuclear fusion. Confinement devices such as Tokamaks and Stellarators \cite{fasoli2016computational, garabedian2003computational, grandgirard2013gyrokinetic, spitzer1958stellarator, imbert2024introduction} employ sophisticated magnetic field configurations to confine high-temperature plasmas, aiming to achieve the conditions necessary for sustained fusion reactions. Understanding and controlling plasma behavior in these devices is challenging due to the complex interplay of charged particles, self-consistent electromagnetic fields, and externally applied magnetic fields.

In this work, we study the time evolution of plasma density using a Vlasov–Poisson system, which models the dynamics of charged particles under the combined influence of a self-consistent electric field and an externally applied magnetic field. In particular, we focus on the behavior of negatively charged particles, namely electrons, while assuming the ions to form a stationary neutral background.
	Under these assumptions, the equations that govern the evolution of the electron density $f(t,\xx,\vvel)$ take the form
	\begin{equation}\label{eq:vlasov_poisson_intro}
		\begin{split}
			&\partial_t f(t,\xx,\vvel)+ \vvel \cdot \nabla_\xx f(t,\xx,\vvel) + \big( \EE(t,\xx) +\vvel \times \BB_{ext}(t,\xx) \big) \cdot \nabla_\vvel f(t,\xx,\vvel) = 0,\\\vspace{0.2cm}
			&-\Delta_\xx \phi(t,\xx) = \rho(t,\xx) - \rho_i(t,\xx), \qquad \EE(t,\xx) = -\nabla_\xx \phi(t,\xx),
		\end{split}
	\end{equation}with $(\xx,\vvel) \in \RR^{2d}$ being the position and velocities, and where the electrons density is defined by
	\begin{equation}\label{eq:charge_density_intro}
		\rho(t,\xx) = \int_{\mathbb{R}^{d_v}} f(t,\xx,\vvel)\,d\vvel,
	\end{equation}
	while $\rho_i$ represents the ion density. 
	In this formulation, $ \EE(t,\xx) $, is the self-consistent electric field obtained from the Poisson equation as the gradient of the electric potential $\phi(t,\xx)$, while $ \BB_{ext}(t,\xx) $ is an external magnetic field. 

The accurate simulation of such systems is notoriously demanding. Plasma dynamics combine turbulence, strong non-linearities, and rapid temporal oscillations with the additional complexity induced by external and self-consistent fields. The resulting multi-scale nature of the problem requires numerical methods that are at the same time stable, accurate, and capable of resolving the problem across different regimes \cite{crouseilles2004numerical,degond2013ap}. To address these challenges, structure-preserving schemes have been developed, providing a framework to capture the essential physical features of magnetized plasma evolution while keeping computational costs under control \cite{belaouar2009asymptotically,cheng2014discontinuous,crouseilles2016multiscale,Chac2023}.
Among the different numerical schemes, we recall semi-Lagrangian schemes and finite volumes ones \cite{coughlin2022efficient, crouseilles2004highorder, crouseilles2010conservative, degond2013ap, degond2017ap, russo2009semilagrangian,sonnendrucker1999semilagrangian,yang2014hermite}. In this work, we will focus on Particle-In-Cells schemes  \cite{deluzet2023parallelization, chacon2016pic,taitano2013implicit,dimarco2010dsmc,filbet2016pic,filbet2017ap,gu2022hamiltonian,hairer2018boris}. Specifically, we approximate the plasma density by a finite number of particles whose trajectories are given by the characteristic curves originating from the Vlasov equation. 

A fundamental aspect in the design of numerical methods for plasma control is the presence of an external magnetic field, which in this context serves as the control variable required to drive the plasma toward the desired configuration. 
More specifically, the control problem we aim to study is given by 
	\begin{equation}\label{eq:control_pb2}
		\min_{\BB_{ext}\in \mathcal{B}_{adm}}\mathcal{J}(\BB_{ext}; f^0;f) ,\qquad  \textrm{subject to }\quad \eqref{eq:vlasov_poisson_intro},
	\end{equation} 
	where $\mathcal{B}_{adm}$ is a set of admissible controls, $f^0$ the initial datum, and $\mathcal{J}(\cdot)$ is a functional which reads as follows
	\begin{equation}\label{eq:control_pb}
		\begin{aligned}
			\mathcal{J}(\BB_{ext};f^0;f)= & \int_{0}^{t_f} \left(  \mathcal{D}(f,\psi)(t)  +  \frac{\gamma}{2}  \int_{\RR^{2d}}  |\BB_{ext}(t,\xx)|^2  f(t,\xx,\vvel) d\xx d\vvel \right)dt
		\end{aligned}
	\end{equation}
	where 
	\begin{equation}\label{eq:D_continuos} 
		\mathcal{D}(f,\psi)(t) = \int_{\RR^{2d}}    \psi(\xx,\vvel) \left( f(t,\xx,\vvel) - \hat{f}(t,\xx,\vvel)  \right) d\xx d\vvel ,
	\end{equation}
	is a running cost function, $\hat{f}$ is a given target distribution discussed next, $\psi(\xx,\vvel)$ is a function of the state variables, $t_f$ is a final prescribed time, and  $\gamma$ a weight penalizing the magnitude of the control given by the external magnetic field $\BB_{ext}$. The scope of the functional in \eqref{eq:control_pb} is to force the distribution function or its moments towards desired values through the choice of the function $\psi(\xx,\vvel)$. Specifically, we aim at reducing the thermal energy and/or the mass at the boundaries of the domain, where excess energy or particle accumulation can create strong gradients or localized instabilities.  
Many efforts have been made in this direction, studying the problem from a theoretical viewpoint as in \cite{bartsch2024controlling, knopf2020optimal, caprino2012magnetic, han-kwan2010tokamak, weber2021optimal}, and from a numerical perspective as in \cite{einkemmer2024vlasov, einkemmer2025control}.
The aim of the present paper is to provide an overview of the results introduced in \cite{albi2025instantaneous,albi2025robust}, reformulating the various control strategies in a two-dimensional framework using polar coordinates. While Cartesian formulations remain widely used, they often fail to exploit the inherent symmetries of confinement devices. In contrast, polar coordinates provide a more natural description of the plasma geometry, especially in Tokamak-like configurations. This choice not only simplifies the representation of magnetic field lines and flux surfaces but also enhances the efficiency and accuracy of numerical discretization.
The idea is to construct a piecewise constant control in space by introducing a computational grid composed of cells. In the approach introduced in \cite{albi2025instantaneous}, the control is directly derived as a piecewise constant function within each cell. This strategy is straightforward to implement and ensures that the control can be incorporated efficiently into the numerical scheme at the cell level. In the second approach introduced in \cite{albi2025robust}, the control is initially computed at the level of individual particles, taking into account the local state of the system. The particle-level information is then interpolated to define an average magnetic field within each cell. This approach allows for a finer resolution of the control, as it captures sub-cell variations in the plasma dynamics before producing a cell-averaged control. Consequently, it can provide a more accurate representation of the magnetic field needed to steer the plasma toward the desired configuration, particularly in regions where strong gradients or localized effects are present. Numerical experiments will show a comparison between the two different approaches, proving the efficiency of both strategies.

The rest of the paper is organized as follows. In Section 
\ref{sec:problem_setting} we formulate the problem in polar coordinates and we introduce the first order semi-implicit Particle-In-Cell scheme that we will employ for simulations. In Section \ref{sec:control} we formulate the two control problems in polar coordinates. In Section \ref{sec:numerical_exp} we conduct different numerical experiments to demonstrate the accuracy and efficiency of the proposed strategies. Finally, in Section \ref{sec:conclusion} we highlight possible future research lines. 
\section{Problem setting and numerical methods}\label{sec:problem_setting} 
In this work, we adopt a simplified framework by considering a two-dimensional phase space. While positions and velocities were initially considered in $\RR^d$, we now focus on a two-dimensional phase space and investigate the plasma dynamics within a circular spatial domain. Specifically, particles are characterized by their position $\xx\in\Omega_x\subset \RR^2$ and velocity  $\vvel\in\Omega_v\subset \RR^2$.
To better capture the geometry of the setting, we express the dynamics in polar coordinates, assigning to each particle a radial coordinate  and an angular coordinate $(r,\theta)\in \Omega_r$, such that
\begin{equation}\label{eq:polar_coords}
	x = r \cos(\theta), \qquad y = r \sin(\theta),
\end{equation}
and the corresponding velocity components are given by
\begin{equation}\label{eq:polar_coords_velocity}
	\begin{split}
		v_x &= v_r \cos(\theta) - v_\theta \sin(\theta),\\
		v_y &= v_r \sin(\theta) + v_\theta \cos(\theta),
	\end{split}
\end{equation}
where $v_r = \dot{r}$ denotes the radial velocity and $v_\theta = r \dot{\theta}$ the angular velocity, with $(v_r,v_\theta)\in \Omega_v$.
We focus on the long time behaviour of the plasma in the orthogonal plane to the external magnetic field which is supposed to be given by
\begin{equation}\label{eq:magnetic_field}
	\BB(t,\rr) = \left( 0,0, B (t,\rr)\right),
\end{equation} 
where $B(t,\rr) >0$. Here and in the following we will use the simplified notation, denoting $\BB_{ext}(\cdot)$ as $\BB(\cdot)$.
The two-dimensional Vlasov-Poisson  equations in polar coordinates reads as 
\begin{equation}\label{eq:vlasov_poisson}
	\begin{split}
		&  \frac{\partial{f(t,\rr,\vvel)}}{\partial{t}} +   \vvel \cdot \nabla_{\rr} f(t,\rr,\vvel) + \left(  \EE(t,\rr) + \frac{v_\theta}{r} \begin{bmatrix}
			v_\theta	\\
			v_r
		\end{bmatrix} + 
		\vvel \times \BB(t,\rr)\right) \cdot \nabla_{\vvel} f(t,\rr,\vvel)   = 0,\\
		& \EE(t,\rr) = -\nabla_{\rr_\perp} \phi(t,\rr),\qquad 	\qquad\nabla_\rr\cdot  \left( \frac{1}{r} \nabla_\rr \phi(t,\rr) \right)   + \frac{1}{r^2} \Delta \phi(t,\rr )  = - \rho(t,\rr),   
	\end{split}
\end{equation}
where $\rr = (r,\theta)$ and $\vvel= (v_r,v_\theta)$. 
In \eqref{eq:vlasov_poisson}, $f(t,\rr,\vvel)$ represents the plasma density which is supposed to lie on the $r-\theta$ plane.

The idea is to control the plasma dynamics by means of an external magnetic field, so as to drive the system toward a desired configuration. To model a more realistic scenario, we assume that the control is piecewise constant in space. Specifically, we introduce a partition of the spatial domain into $N_c$ cells and assign to each cell  $C_k\subset \Omega_r$ a constant control value $B_k$. The cells satisfy $\bigcup_{k=1}^{N} C_k=\Omega_r$, $C_k\cap C_\ell=\emptyset$ for all $k\neq \ell$, and $k, \ell\in\{1,\ldots,N_c\}$, with in the following $\Omega_k = C_k \times \Omega_v$ the state space of the single cell.  We emphasize that the choice of  $N_c$ is not dictated by numerical discretization, but rather reflects the physical constraints of the problem. Indeed, from a practical perspective, one cannot expect the coils generating the external magnetic field to reproduce highly complex or pointwise-varying field structures.

We start from the following continuous control problem
\begin{equation}\label{eq:control_pb_continuos}
	\min_{\BB \in \mathcal{B}_{adm}}    \mathcal{J}(\BB;f^0,f) ,\qquad 
	\textrm{s.t. }~\eqref{eq:vlasov_poisson} \ \text{is satisfied},
\end{equation}
where 
\begin{equation}\label{eq:functional_continuos}
	\begin{split}
		&\mathcal{J}(\BB;f^0,f)  = \int_{0}^{t_f}  \sum_{\ell \in \{\text{r},\text{v}\}}\mathcal{D}(f,\psi_\ell)(t)   dt +\cr
		&\qquad + \frac{\gamma}{2} \int_{0}^{t_f}  \int_{\Omega} |\BB(t,\rr) |^2 f(t,\rr,\vvel)d\rr d\vvel  dt,
	\end{split}
\end{equation}
with $\Omega = \Omega_r \times \Omega_v$, and 
where %$\gamma>0$ is a penalization term, $\mathcal{B}_{adm}$ is a set of admissible controls,  $\mathcal{P}[\cdot]$ is a suitable statistical operator taking into account the presence of the uncertainties, and 
$\mathcal{D}(\cdot)$ aims at enforcing a specific configuration of the distribution function and of its moments, that is 
\begin{equation}\label{eq:D}
	\begin{split}
		\mathcal{D}(f,\psi_\ell)(t) = &\frac{\alpha_\ell}{2} | m(f,\psi_\ell)(t) - \hat{\psi}_\ell |^2 + \frac{\beta_\ell}{2} m_\sigma(f,\psi_\ell)(t),
	\end{split}
\end{equation}
with $\alpha_\ell,\beta_\ell \geq 0$ weighting parameters for $\ell = \{\text{x},\text{v}\}$, and
\begin{equation}\label{eq:m_sigma} 
	\begin{split}
		&m(f,\psi_\ell)(t) = \int_{\Omega} \psi_\ell(\rr,\vvel) f(t,\rr,\vvel) d\rr d\vvel,\\
		& m_{\sigma}(f,\psi_\ell) (t) = \int_{\Omega} | \psi_\ell(\rr,\vvel)-m(f,\psi_\ell)(t)|^2 f(t,\rr,\vvel) d\rr d\vvel.
	\end{split}
\end{equation}
The first equation in \eqref{eq:m_sigma} is designed to enforce a desired configuration of the distribution function corresponding to a prescribed target $\hat{\psi}_\ell$. The second equation minimizes the distance from the average value $m(f,\psi_\ell)$ both in space as well as in velocity space around the target. 
In this simplified framework, the goal is to prevent the mass from reaching the upper and lower boundaries of the domain, aiming instead for a particle configuration concentrated near the center with minimal variance.  

In the following Sections, we will present an overview on the different control problems introduced in \cite{albi2025instantaneous,albi2025robust}, reformulated in polar coordinates.  

\subsection{PIC scheme in polar coordinates }
At the numerical level, the density function $f$ is approximated by a set of $N$ particles with state $\rr_m$, $\vvel_m$, for $m=1,\ldots,N$ as follows
	\begin{equation}\label{eq:approx_density}
		f^N(t,\rr,\vvel) = \frac{1}{N}\sum_{m=1}^N \delta(\rr-\rr_m(t))\delta(\vvel-\vvel_m(t)),
\end{equation}
where $\delta(\cdot)$ is the Dirac-delta function.
The trajectory of each particle $m=1,\ldots,N$ is computed from the characteristic curves originated by the Vlasov equation in \eqref{eq:vlasov_poisson},  
\begin{equation}\label{eq:char_curves_polar} 
	\begin{split}
		& \frac{dr_m(t)}{dt} = v_{r_m}(t), \qquad r_m(0) = r^0_m,  \\
		& \frac{d\theta_m(t)}{dt} = \frac{v_{\theta_m}(t)}{r_m(t)}, \qquad v_{\theta_m}(0) = v_{\theta_m}^0,  \\
		& \frac{dv_{r_m}(t)}{dt} = E_{r_m}(t) + v_{\theta_m}(t) \left(  B_m(t) +\frac{v_{\theta_m}(t)}{r_m(t)} \right)  , \qquad v_{r_m}(0) = v_{r_m}^0,  \\
		& \frac{dv_{\theta_m}(t)}{dt} = E_{\theta_m}(t) - v_{r_m}(t) \left(  B_m(t) +\frac{v_{\theta_m}(t)}{r_m(t)} \right)  , \qquad v_{\theta_m}(0) = v_{\theta_m}^0,  
	\end{split}
\end{equation}
being $B_m(t) = B(t,\rr_m)$ the magnetic field, and $E_{r_m}(t)$, $E_{\theta_m}(t)$ the radial and angular component of the electric field associated to each particle $m=1,\ldots,N$. The electric field is computed over a spatial grid by means of a finite difference scheme as the solution to the Poisson equation in \eqref{eq:vlasov_poisson}. Then the value of the electric field associated to each particle is taken as the value at the center of the cell in the spatial discretization grid containing it. At each time step,
the approximated density, needed to solve the Poisson equation, is reconstructed over the
spatial cells from the updated particle state as in \eqref{eq:approx_density}. We remark that the spatial discretization grid and the fictitious cells on which the control takes constant values are different. One is dictated by the numerical discretization while the other by the physics of the problem, as specified in the following section.  
We then consider a time interval $[0,T]$ divided into $N_t$ time steps of size $h$ and we define $t^n = n h$. To discretize \eqref{eq:char_curves_polar}  we rely on a first order semi-implicit scheme
\begin{equation}\label{eq:PIC_scheme}
	\begin{split}
		&  r_m^{n+1} = r_m^n + h v_{r_m}^{n+1},  \\
		&  \theta_m^{n+1} =\theta_m^{n} + h \frac{v_{\theta_m}^{n+1}}{r_m^{n+1}},  \\
		&  v_{r_m}^{n+1} =  v_{r_m}^{n} + h \left[  E^n_{r_m} + v_{\theta_m}^{n+1} \left(  B_m^n +\frac{v_{\theta_m}^n}{r_m^n} \right) \right], \\
		&  v_{\theta_m}^{n+1} = v_{\theta_m}^{n}+ h\left[  E_{\theta_m}^n - v_{r_m}^{n+1} \left(  B_m^n +\frac{v_{\theta_m}^n}{r_m^n} \right)  \right] ,  
	\end{split}
\end{equation}
where $B_m$ is set either to be constant or it is computed as the solution of one of the control problems introduced in the following section.

\section{Control strategies}\label{sec:control} 
In this section, we present an overview on the control strategies introduced in \cite{albi2025instantaneous,albi2025robust} based on the derivation of the magnetic field $B$ as a control taking piecewise constant values on the physical space. We show the derivation in polar coordinates. 

We consider a fictitious space discretization grid with $N_c$ cells $C_k$ of size $\Delta^c  r \times \Delta^c \theta$ with $\Delta^c r\gg\Delta r$, and $\Delta^c \theta \gg\Delta \theta$, being $\Delta r$, $\Delta \theta$ the size of the space discretization grid on which the density is reconstructed. %(see Figure \ref{fig:griglia}), 
Then, we show how to derive a feedback control $B_k$, $k=1,\ldots,N_c$ taking constant values in each cell $C_k$ at every instant of time $t^n$, and based on a one-step prediction of the dynamics.

Strategy one consists in directly computing the piecewise constant control as the solution of the underlying optimization problem. In contrast, strategy two first determines a pointwise control, i.e., a control value associated with each particle, and subsequently projects it onto the piecewise constant structure through interpolation.

\subsection{Strategy one} \label{sec:strategy1}
We first formulate the problem at the continuous level and over a finite time horizon $[0,t_f]$ as follows
\begin{equation}\label{eq:continuos_pb_strategy1}
	\min_{B\in \mathcal{B}_{adm}}  \sum_{k=1}^{N_c} \mathcal{J}_k(B_k;f_k,f^0_k),\qquad 
	\textrm{s.t.}~\eqref{eq:vlasov_poisson},
\end{equation}
where  $f_k= f_k (t,\rr,\vvel)$ corresponds to the normalized particle density restricted to a single cell $C_k$ for each $k = 1,\ldots, N_c$
\[
f_k (t,\rr,\vvel) = \frac{f(t,\rr,\vvel)}{\rho_k(t)},\qquad \rho_k(t) = \int_{\Omega_k}f(t,\rr,\vvel)\,d\rr \, d\vvel, 
\]
with $\rho_k(t)>0$ the total cell density and with $B=(B_1,\ldots,B_{N_c})$ now representing the vector of $z$ components of $\BB(t,\rr)$   within each cell $C_k$, $\mathcal{B}_{adm}$ the set of admissible controls such that
$\mathcal{B}_{adm} = \{B_k  | B_k \in[-M,M], M>0,\, k = 1,\ldots, N_c\}$. For each $k=1,\ldots,N_c$, the cost functional in \eqref{eq:continuos_pb_strategy1} is defined as follows
\begin{equation}\label{eq:J_strategy1}
	\begin{split}
		\mathcal{J}_k(B_k; f_k,f_k^0) = & \int_{0}^{t_f} \left( \sum_{\ell\in\{\textrm{r},\textrm{v}\}} \mathcal{D}(f_k)(t,\psi_\ell)  + \right. \\ &  \qquad + \left. \frac{\gamma}{2}  \int_{\Omega_k} | B(t,\rr)|^2 f_k(t,\rr,\vvel) d\rr d\vvel \right)\, dt,
	\end{split}
\end{equation}
where $\gamma>0$ is a penalization
term,  $\mathcal{D}(\cdot)$ aims at enforcing a specific configuration in the distribution function, and $\Omega_k = C_k \times \Omega_v$, being $\Omega_v$ the velocity domain.  

We consider a short time horizon of
length $h > 0$ and formulate a time discretize optimal control problem through the
functional $\mathcal{J}_k$ restricted to the interval $[t, t + h]$, as follows
\begin{equation}\label{eq:min_prob_cell_strategy1}
	\min_{B_k\in \mathcal{B}_{adm}} \mathcal{J}^{N,h}_{k}(B_k; f_k^N,f_k^{N,0}),
\end{equation}
subject to a semi-implicit in time discretized Vlasov dynamics, fully explicit for the velocity terms
\begin{equation}\label{eq:explicit_dinamics}
	\begin{split}
		&r_m^{n+1}  = r_m^n  + h v_{r_m}^{n+1} ,\\
		&\theta_m^{n+1} = \theta_m^n + h \frac{v_{\theta_m}^{n+1}}{r_m^{n+1}},\\
		&v_{r_m}^{n+1}  = v_{r_m}^n + h  \frac{(v_{\theta_m}^n)^2}{r_m^n} + v_{\theta_m}^n B_k  + h E_{r_m}^n,\\
		& v_{\theta_m}^{n+1}  = v_{\theta_m}^n + h  \frac{v_{\theta_m}^n v_r^n}{r_m^n} - v_{r_m}^n B_k + h E_{\theta_m}^n.
	\end{split}
\end{equation}
Using the rectangle rule for approximating the integral in time, the functional in \eqref{eq:min_prob_cell_strategy1} reads as follows
\begin{equation}\label{eq:discr_J_1_strategy1}
	\mathcal J_k^{N,h}(B_k;f_k^N,f_k^{N,0}) =h  \left( \sum_{\ell\in\{\textrm{r},\textrm{v}\}} \mathcal{D}(f_k^{N,n+1},\psi_\ell)(t)+ \frac{\gamma}{2} |B^{n+1}_k|^2 \right).
\end{equation}

Thus, by setting $\psi_\textrm{r} = r^{n+1}$, $\psi_{\textrm{v}} = v_r^{n+1}$, $\hat{\psi}_{\textrm{r},k} = \hat{r}_k$ and $\hat{\psi}_{\textrm{v},k} = \hat{v}_{r_k}$ target states,  and by direct computation over the empirical densities, we can rewrite the functional in \eqref{eq:discr_J_1_strategy1} as
\begin{equation}\label{eq:discr_J_strategy1}
	\begin{split}
		\mathcal J_k^{N,h}(B_k)= & \frac{h\alpha_\emph{v}}{2}  \vert \bar{v}_{r,k}^{n+1}  -\hat{v}_{r_k}\vert^2 + \frac{h\beta_\emph{v}}{2 N_k} \sum_{i\in C_k} \vert v_{r_i}^{n+1} -\bar{v}_{r,k}^{n}  \vert^2 +  \\
		&  + \frac{h\alpha_\emph{r}}{2}  \vert \bar{r}_{k}^{n+1}  -\hat{y}_k\vert^2 + \frac{h\beta_\emph{r}}{2 N_k} \sum_{i\in C_k} \vert r_{i}^{n+1} -\bar{r}_k^{n} \vert^2+  \frac{h\gamma}{2}  \vert B^{n+1}_k \vert^2   ,
	\end{split}
\end{equation}
with 
\begin{equation}\label{eq:mean_quantities}
	\begin{split} 	
		\bar{r}_{k}   = \frac{1}{N_{k}} \sum_{j\in C_k} r_{j}  ,\qquad
		\bar{v}_{r,k}   = \frac{1}{N_{k}} \sum_{j\in C_k} v_{r_j},
	\end{split}
\end{equation}
denoting the mean position and velocity over cell $C_k$, $k=1,\ldots,N_c$.
We reformulate now the result proved in \cite{albi2025instantaneous} in polar coordinates. 
\begin{proposition}\label{prop:istctrl_1}
	Assume the parameters to scale as 
	\begin{equation}\label{eq:scaling}
		\alpha_\emph{r} \rightarrow \frac{\alpha_\emph{r}}{h}, \qquad \beta_\emph{r} \rightarrow \frac{\beta_\emph{r}}{h}, \qquad \gamma \rightarrow \gamma h,  
	\end{equation}
	then the feedback control at cell $C_k$ associated to \eqref{eq:discr_J_strategy1} reads as follows
	\begin{equation}\label{eq:Bk_strategy1_h}
		B_k = \mathbb{P}_{[-M,M]}\left( \frac{ \mathcal{R}^{N,n}_{\emph{v},k}   + \mathcal{R}^{N,n}_{\emph{r},k}   }{\gamma + \mathcal{S}^{N,n}_{\emph{v},k}  + \mathcal{S}^{N,n}_{\emph{r} ,k}  }\right), 
	\end{equation}
	where $\gamma>0$, 
	\begin{equation}\label{eq:terms_in_B}
		\begin{split}
			\mathcal{R}^{N,n}_{\emph{v},k} = & \alpha_\emph{v}  \bar{v}_{\theta,k}^n \left(\bar{v}_{r,k}^n + h \frac{(\bar{v}_{\theta,k}^n)^2}{\bar{r}_k^n} + h \bar{E}_{r,k}^n - \hat{v}_{r_k}\right) + \\
			& \qquad + \frac{\beta_\emph{v}}{N_k} \sum_{i\in C_k} \left[  v_{\theta_i}^n \left(v_{r_i}^n + h \frac{(v_{\theta_i}^n)^2}{r_i^n} + h E_{r_i}^n - \bar{v}_{r_k}^n\right)\right]  ,\\
			\mathcal{R}^{N,n}_{\emph{r},k}    =  & \alpha_\emph{r} \bar{v}_{\theta,k}^n \left( \bar{r}_{k}^n + h \bar{v}_{r,k}^n + h^2 \frac{(\bar{v}_{\theta,k}^n)^2}{\bar{r}_k^n} + h^2 \bar{E}_{r,k}^n - \hat{r}_k  \right) + \\
			&  \frac{\beta_\emph{r}}{N_k} \sum_{i\in C_k} \left[v_{\theta_i}^n \left( r_{i}^n + h v_{r_i}^n + h^2 \frac{(v_{\theta_i}^n)^2}{r_i^n} + h^2 E_{r_i}^n - \bar{r}_k^n  \right)\right],\\
			\mathcal{S}^{N,n}_{\emph{v},k}  = & h  (\alpha_\emph{v} +\beta_\emph{v}) (\bar{v}_{\theta,k}^n)^2,\\
			\mathcal{S}^{N,n}_{\emph{r},k}  = &  h^2 (\alpha_\emph{r} +\beta_\emph{r}) (\bar{v}_{\theta,k}^n)^2,
		\end{split}
	\end{equation}
	and with $\mathbb{P}_{[-M,M]}(\cdot)$ denoting the projection over the interval $[-M,M]$. 
	In the limit $h\to 0$ the control at the continuous level reads, 
	\begin{equation}\label{eq:Bk_strategy1}
		B_k(t) = \mathbb{P}_{[-M,M]} \left( \frac{1}{\gamma} \left(	\mathcal{R}^N_{\emph{v},k} (t) + 	\mathcal{R}^{N}_{\emph{r},k} (t)  \right) \right),
	\end{equation}
	with 
	\begin{equation*}
		\begin{split}
			&	\mathcal{R}^{N}_{\emph{v},k} (t) = \alpha_\emph{v} \bar{v}_{\theta,k}(t) (\bar{v}_{r,k}(t) - \hat{v}_{r_k}) + \frac{\beta_\emph{v}}{N_k} \sum_{i\in C_k} \left[  v_{\theta_i}(t) (v_{r_i}(t) - \bar{v}_{r_k}(t))  \right] ,\\
			&	\mathcal{R}^{N}_{\emph{r},k} (t) = \alpha_\emph{r} \bar{v}_{\theta,k}(t) (\bar{r}_{r,k}(t) - \hat{r}_{k}) + \frac{\beta_\emph{r}}{N_k} \sum_{i\in C_k} \left[  v_{\theta_i}(t) (r_{i}(t) - \bar{r}_{k}(t))  \right].
		\end{split}
	\end{equation*}
\end{proposition}
The proof follows the approach outlined in \cite{albi2025instantaneous}, where we show how to solve the corresponding optimality system using an augmented Lagrangian method.

\subsection{Strategy two} \label{sec:strategy2}
We consider the minimization problem in \eqref{eq:control_pb_continuos} and as before we formulate a time discretize optimal control problem restricted to the interval $[t,t+h]$ as follows 
\begin{equation}\label{eq:min_prob_cell_strategy2}
	\min_{B_m\in \mathcal{B}_{adm}} \mathcal{J}_m^{N,h}(B_m; f^N,f^{N,0}),
\end{equation}
subject to the semi-implicit in time discretized Vlasov dynamics, fully explicit for the velocity terms given by \eqref{eq:explicit_dinamics} with $B_k$ substituted by $B_m$. In \eqref{eq:min_prob_cell_strategy2} $\mathcal{B}_{adm}$ is the set of admissible controls such that $\mathcal{B}_{adm} = \{B_m  | B_m \in[-M,M], M>0,\, m = 1,\ldots, N\}$, and the functional reads as 
\begin{equation}\label{eq:J_strategy2}
	\begin{split}
		\mathcal{J}_m(B_m; f^N,f^{N,0}) = & \int_{0}^{t_f} \left( \sum_{\ell\in\{\textrm{r},\textrm{v}\}} \mathcal{D}(f^N)(t,\psi_\ell)  + \right. \\ &  \qquad + \left. \frac{\gamma}{2}  \int_{\Omega} | B(t,\rr)|^2 f^N(t,\rr,\vvel) d\rr d\vvel \right)\, dt,
	\end{split}
\end{equation}
Then, we consider a short time horizon of length $h>0$ and we formulate the problem in the interval $[t,t+h]$ as before. Applying a semi-implicit rectangle rule to approximate the time integrals in \eqref{eq:J_strategy2}, we obtain 
\begin{equation}\label{eq:J_discretized_h} 
	\mathcal{J}_m^{N,h}(B_m;f^{N},f^{N,0}) = h  \left(   \sum_{\ell \in \{\text{r},\text{v}\}} \mathcal{D}(f^{N,n+1},\psi_\ell)(t)+ \frac{  \gamma}{2} |B_m^{n+1}|^2 \right) .
\end{equation}
By direct computation assuming in \eqref{eq:D} $\psi_\emph{r} = r^{n+1}$, $\psi_\emph{v} = v_r^{n+1}$, $\hat{\psi}_\emph{r} = \hat{r}$,  $\hat{\psi}_\emph{v} = \hat{v}_r$, being $\hat{r}, \hat{v}_r$  certain target functions, we get 
\begin{equation}\label{eq:J_discr_strategy2} 
	\begin{split}
		&\mathcal{J}^{h,N}(B_m) =   h \sum_{m=1}^N \left(  \frac{\alpha_\emph{r}}{2} \Big \vert r_m^{n+1}  -\hat{r} \Big \vert^2 + \frac{\beta_\emph{r}}{2} \Big \vert r_m^{n+1} -\bar{r}^n \Big \vert^2 \right.   + \\
		& \qquad + \left. \frac{\alpha_\emph{v}}{2} \Big \vert v_{r_m}^{n+1} - \hat{v}_r \Big \vert^2 + \frac{\beta_\emph{v}}{2} \Big \vert v_{r_m}^{n+1}  -\bar{v}_r^n\Big \vert^2   
		+ \frac{ \gamma}{2}    \vert B_m^{n+1}\vert^2 \right) ,
	\end{split}
\end{equation}
with $\alpha_\emph{r},\alpha_\emph{v}, \beta_\emph{r},\beta_\emph{v}, \gamma \geq 0$, and where 
\begin{equation*}
	\bar{r}^n  = \frac{1}{N} \sum_{m=1}^N r_m^n , \qquad 	\bar{v}_r^n  = \frac{1}{N} \sum_{m=1}^N v_{r_m}^n.
\end{equation*} 
We now reformulate the result established in \cite{albi2025robust} in polar coordinates.
\begin{proposition}\label{prop:istctrl_2}
	Assume the parameters to scale as 
	\begin{equation}\label{eq:scaling_new}
		\alpha_\emph{r} \rightarrow \frac{\alpha_\emph{r}}{h}, \qquad \beta_\emph{r} \rightarrow \frac{\beta_\emph{r}}{h}, \qquad \gamma \rightarrow \gamma h, 
	\end{equation}
	then the feedback control $B_m$ associated to \eqref{eq:J_discr_strategy2} reads as follows
	\begin{equation}\label{eq:Bk_strategy2_h}
		B_m = \mathbb{P}_{[-M,M]}\left( \frac{ \mathcal{R}^{N,n}_{\emph{r},m}   + \mathcal{R}^{N,n}_{\emph{v},m}  }{\gamma + \mathcal{S}^{N,n}_{\emph{r},m}  + \mathcal{S}^{N,n}_{\emph{v},m}  } \right), 
	\end{equation}
	where for any $m=1,\ldots,N$, $\gamma>0$,	\begin{equation}\label{eq:terms_in_B_strategy2}
		\begin{split}
			\mathcal{R}^{N,n}_{\emph{v},m} = & \alpha_\emph{v}    v_{\theta_m}^n \left(v_{r_m}^n + h \frac{(v_{\theta_m}^n)^2}{r_m^n} + h E_{r_m}^n - \hat{v}_{r_m}\right) + \\
			& \qquad +\beta_\emph{v}   v_{\theta_m}^n \left(v_{r_m}^n + h \frac{(v_{\theta_m}^n)^2}{r_m^n} + h E_{r_m}^n - \bar{v}_{r}^n\right)  ,\\
			\mathcal{R}^{N,n}_{\emph{r},m}    =  & \alpha_\emph{r} v_{\theta_m}^n \left(  r_m^n + h v_{_m}^n + h^2 \frac{(v_{\theta_m}^n)^2}{r_m^n} + h^2 E_{r_m}^n - \hat{r}  \right) + \\
			& \beta_\emph{r}  v_{\theta_m}^n \left( r_{m}^n + h v_{r_m}^n + h^2 \frac{(v_{\theta_m}^n)^2}{r_m^n} + h^2 E_{r_m}^n - \bar{r}^n  \right),\\
			\mathcal{S}^{N,n}_{\emph{v},m}  = & h  (\alpha_\emph{v} +\beta_\emph{v}) ( v_{\theta_m}^n)^2,\\
			\mathcal{S}^{N,n}_{\emph{r},m}  = &  h^2 (\alpha_\emph{r} +\beta_\emph{r}) (v_{\theta_m}^n)^2, 
		\end{split}
	\end{equation}
	and with $\mathbb{P}_{[-M,M]}(\cdot)$ denoting the projection over the interval $[-M,M]$. 
	In the limit $h\to 0$ the control at the continuous level reads, 
	\begin{equation}\label{eq:Bk_strategy2}
		B_k(t) = \mathbb{P}_{[-M,M]} \left( \frac{1}{\gamma} \left(	\mathcal{R}^N_{\emph{v},m} (t) + 	\mathcal{R}^{N}_{\emph{r},m} (t)  \right) \right),
	\end{equation}
	with 
	\begin{equation*}
		\begin{split}
			&	\mathcal{R}^{N}_{\emph{v},m} (t) = \alpha_\emph{v}  v_{\theta_m}(t) (v_{r_m}(t) - \hat{v}_{r}) + \beta_\emph{v} v_{\theta_m}(t) (v_{r_m}(t) - \bar{v}_{r}(t))   ,\\
			&	\mathcal{R}^{N}_{\emph{r},m} (t) = \alpha_\emph{r} v_{\theta_m}(t) (r_{m}(t) - \hat{r}) +\beta_\emph{r} v_{\theta_m}(t) (r_{m}(t) - \bar{r}(t)).
		\end{split}
	\end{equation*}
\end{proposition}
The proof follows the approach described in \cite{albi2025robust}, and is analogous to the one of Proposition \ref{prop:istctrl_1}.
\\\\
Proposition \ref{prop:istctrl_2} establishes a feedback control law based on the pointwise evaluation of the magnetic field at each particle’s position. Such a requirement, however, is not technologically feasible in realistic applications. To relax this constraint, we make use of the previously introduced partition of the domain into cells $C_k$, $k=1,\ldots,N_c$. 
For each $k = 1, \ldots, N_c$, we define the piecewise constant control  $\hat{B}_k$ over the cell $C_k$ by interpolating the pointwise control values $B_m$  within the corresponding cell. For convenience, we denote by  
\begin{equation}\label{eq:Bk_strategy2_interp}
	\BB = \mathcal{I}(\BB_N, \rr_c),
\end{equation}
the vector $\BB = [\hat{B}_1, \ldots, \hat{B}_{N_c}]$ collecting the interpolated controls $\hat{B}_k$. Here, $\mathcal{I} (\cdot)$ denotes a piecewise-constant interpolation operator, $\rr_c$ is the vector of the positions of the centres of the cells $C_k$, and $\BB_N = [B_1, \ldots, B_N]$ is the vector of pointwise controls $B_m$ as defined in equation~\eqref{eq:Bk_strategy2_h}. In this framework, the magnetic field acting on a particle located at position $\rr_m$ is approximated by the constant field value assigned to the cell $C_k$ containing the particle. While this approach increases the complexity of the control problem, it provides a more realistic and practically implementable setting.
\begin{remark}
	Both controls defined in \eqref{eq:Bk_strategy1}–\eqref{eq:Bk_strategy2} are suboptimal, as they are obtained from an explicit discretization of the Vlasov equation \eqref{eq:explicit_dinamics} and subsequently incorporated into the first-order semi-implicit scheme \eqref{eq:PIC_scheme} to steer the particles toward the desired configuration.
\end{remark}
	\begin{remark}
		Note that the instantaneous controls derived in \eqref{eq:Bk_strategy1_h}–\eqref{eq:Bk_strategy2_h}  are of order $\mathcal{O}(h^2)$. To obtain an effective contribution of the control in the limit $h \to 0$, we assume that the parameters scale as in \eqref{eq:scaling}-\eqref{eq:scaling_new} with the time step $h$, so that the leading-order is recovered \cite{albi2017mean}. Additionally, as discussed in Remark 1 of \cite{albi2025instantaneous}, the same instantaneous control laws can be derived directly from the discretized functional over the interval $[t,t+h]$ without explicitly performing any rescaling. This highlights that the scaling is a mathematical strategy to ensure a non-vanishing control, rather than a physical requirement.
\end{remark}
\section{Numerical experiments: Diocotron instability} \label{sec:numerical_exp} 
\begin{figure}[h!]
	\centering
	\includegraphics[width=0.458\linewidth]{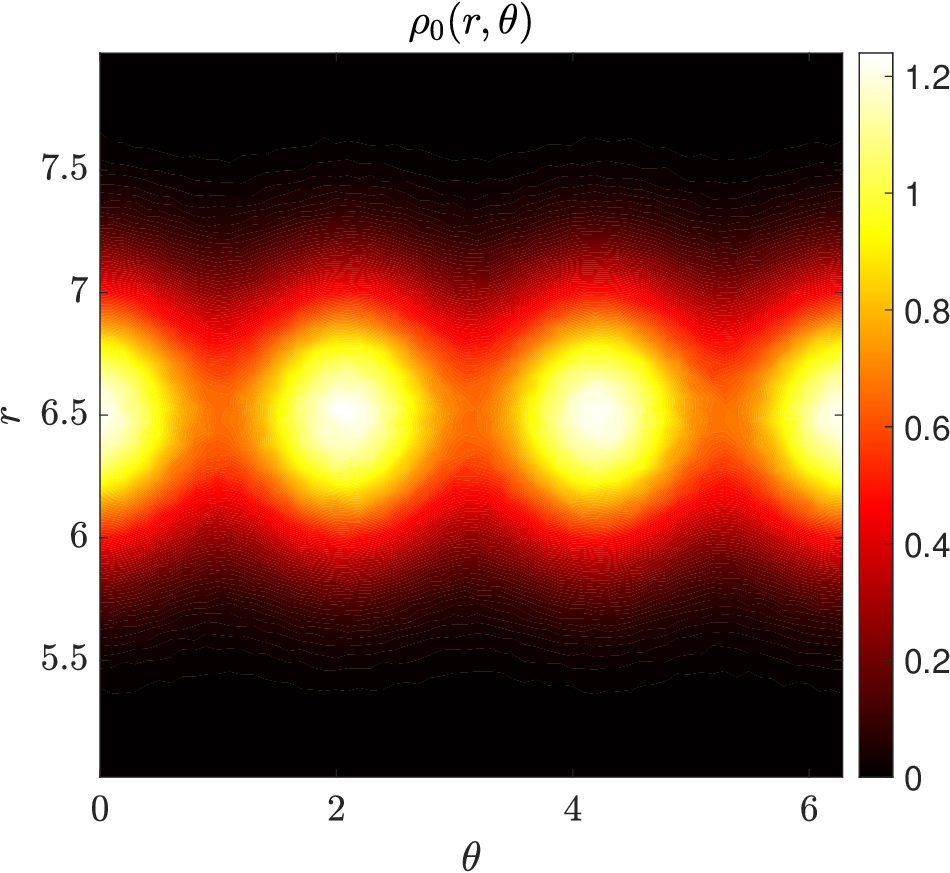}
	\includegraphics[width=0.44\linewidth]{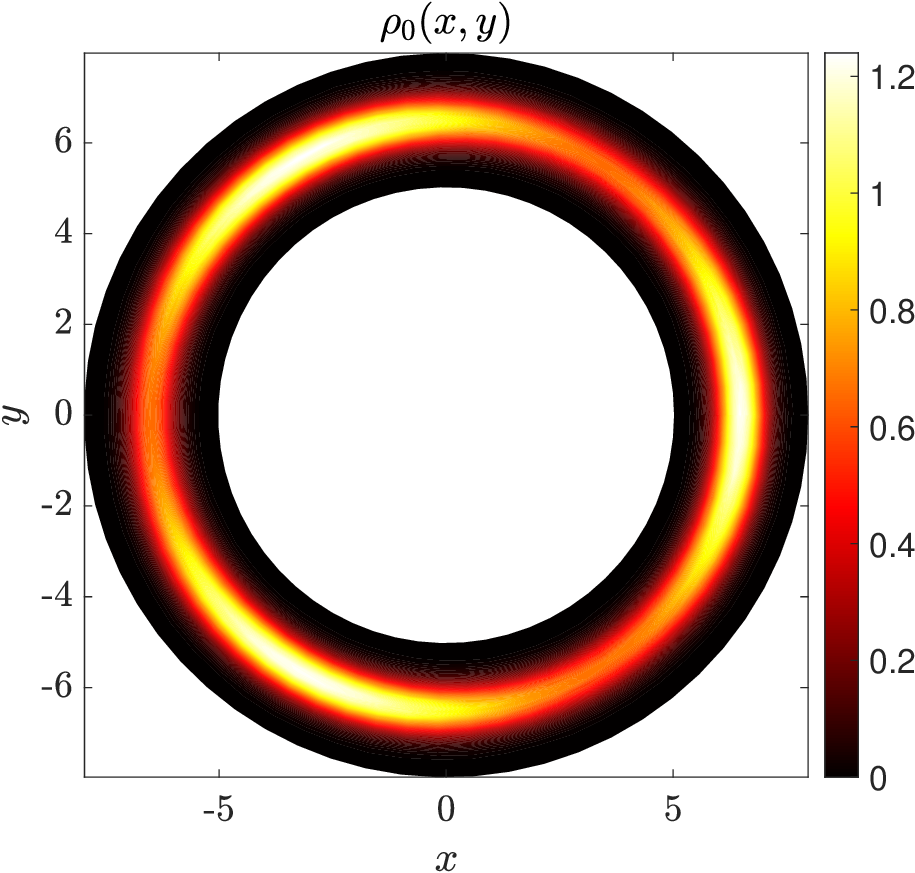}
	\caption{Diocotron instability. Initial density in polar coordinates (on the left) and in cartesian coordinates (on the right). }
	\label{fig:ic}
\end{figure}

In this section, we present a series of numerical experiments to illustrate the effectiveness of the instantaneous control strategies introduced above in steering the plasma toward desired configurations. Our attention is devoted in particular to the Diocotron instability test \cite{crouseilles2013semi,valentini2005numerical}, with the following choice of initial condition 
\begin{equation}\label{eq:diocotron}
	f_0(\rr,\vvel)	= \rho_0(\rr)\exp{\left(- \frac{\vvel^2}{2}\right)},
\end{equation} 
with 
\begin{equation}\label{eq:rho0}
	\rho_0(\rr) = \frac{1}{2\pi}\left( 1+\alpha \cos\left( k \theta  \right) \right) \exp{\left( -4\left( r-6.5\right)^2\right)},
\end{equation}
where $\alpha = 0.3$, $k = 3$, $\rr = (r,\theta) \in [5,8]\times [0,2\pi]$, $\vvel\in \mathbb{R}^{2}$. We set Dirichlet boundary conditions in $r$ and periodic boundary condition in $\theta$. Figure \ref{fig:ic} shows the plot of the initial density function \eqref{eq:rho0} in the $r-\theta$ plane (on the left) and in the $x-y$ plane (on the right). 

We begin with a validation experiment to confirm that the PIC scheme defined in \eqref{eq:PIC_scheme} achieves first-order accuracy. Next, we consider the uncontrolled dynamics, assuming that the system evolves under a magnetic field constant in both space and time. Finally, we compare these results with those obtained by replacing the constant magnetic field with the one derived in Proposition \ref{prop:istctrl_1}–\ref{prop:istctrl_2}. To asses the effectivenesses of the proposed control strategy we measure the thermal energy at the boundaries as
	\begin{equation}\label{eq:energy}
		\mathcal{E}_b^n  = \frac{1}{2 N \vert \Omega_b\vert} \sum_{\mathcal{C}_j \in \Omega_b} \sum_{m=1}^N \left\vert \vvel_m^n  - U_b^n  \right\vert^2 \chi(\rr_m^n \in \mathcal{C}_j),
	\end{equation}
	where $\chi(\cdot)$ denotes the indicator function and  
	\begin{equation}\label{eq:mass_boundary}
		\begin{split}
			U_b^n &= \frac{1}{N \vert \Omega_b\vert \rho_b^n } \sum_{\mathcal{C}_j \in \Omega_b} \sum_{m=1}^N \vvel_m^n \chi(\rr_m^n \in \mathcal{C}_j), \\
			\rho_b^n  &= \frac{1}{N \vert \Omega_b\vert} \sum_{\mathcal{C}_j \in \Omega_b} \sum_{m=1}^N \chi(\rr_m^n \in \mathcal{C}_j),
		\end{split}
	\end{equation}
	and   $ \mathcal{C}_j \in \Omega_b $, with $ \Omega_b = [5 ,5+\Delta r] \cup [8-\Delta r, 8] $. This choice corresponds to a region near the  boundaries with width equal to one cell size $\Delta r$. We also introduce the quantity
	\begin{equation}\label{eq:electric_energy} 
		\mathcal{E}_{el}^n = \frac{1}{2 N \vert \Omega_c\vert}  \sum_{\mathcal{C}_j \in \Omega_c}\sum_{m=1}^N \vert \EE_m^n\vert^2 \chi(\rr_m^n \in \mathcal{C}_j)
	\end{equation}
	which represents the energy associated with the self-induced electric field measured at the center of the domain, namely $\Omega_c = [6.5-\Delta r, 6.5+\Delta r]$.

The simulations are performed using $ N = 10^7 $ particles and a grid of $ m_r \times m_\theta $ cells, with $ m_r = m_\theta = 64 $, for the reconstruction of macroscopic quantities. The final time is set to $ t_f = 250 $, with a time step $ h = 0.5 $.
\subsection{Validation test} \label{sec:validation} 
\begin{figure}[h!]
	\centering
	\includegraphics[width=0.458\linewidth]{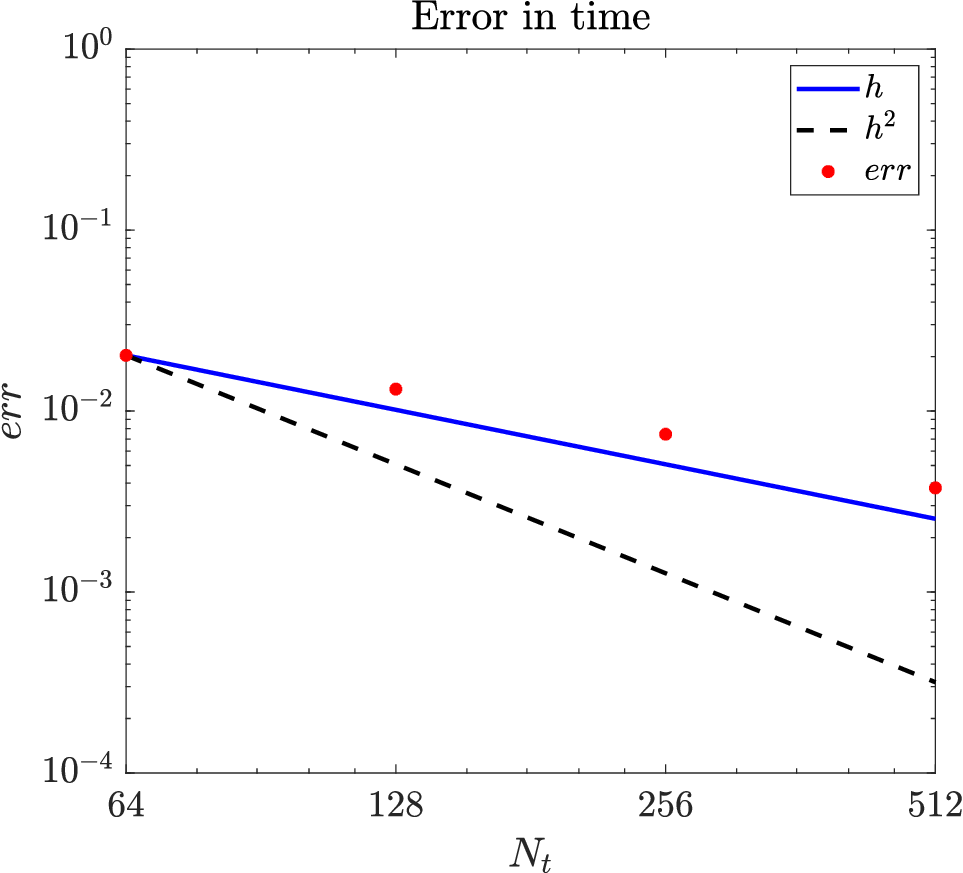}
	\caption{Diocotron instability. Error in time defined as in \eqref{eq:error_time}.  }
	\label{fig:error_time}
\end{figure}

We begin by assessing the temporal accuracy of the scheme introduced in \eqref{eq:PIC_scheme}. To this end, we compute the error in time by comparing the numerical solution against a highly resolved reference solution. Denoting by $U = [\rr_m,\vvel_m]$ the state vector containing the position and velocity of each particle, we define the error at final time $t_f$ as
\begin{equation}\label{eq:error_time}
	err  = \Vert U_{rif} - U \Vert_{\infty},
\end{equation}
where $U_{rif}$ represents the reference solution at time $t_f = 250$ computed with $N_t =2^{12}$ time steps, and $U$ is the corresponding solution computed with coarser resolutions, namely  $N_t= 2^6,\ldots, 2^9$.  Figure \ref{fig:error_time} displays the error computed as in \eqref{eq:error_time} as the number of time steps $N_t$ increases. As expected, the error decreases proportionally to the time step size $h$, confirming that the proposed scheme achieves first-order accuracy in time.

\subsection{Test without control}\label{sec:test_nocontrol}
\begin{figure}[h!]
	\centering
	\includegraphics[width=0.328\linewidth]{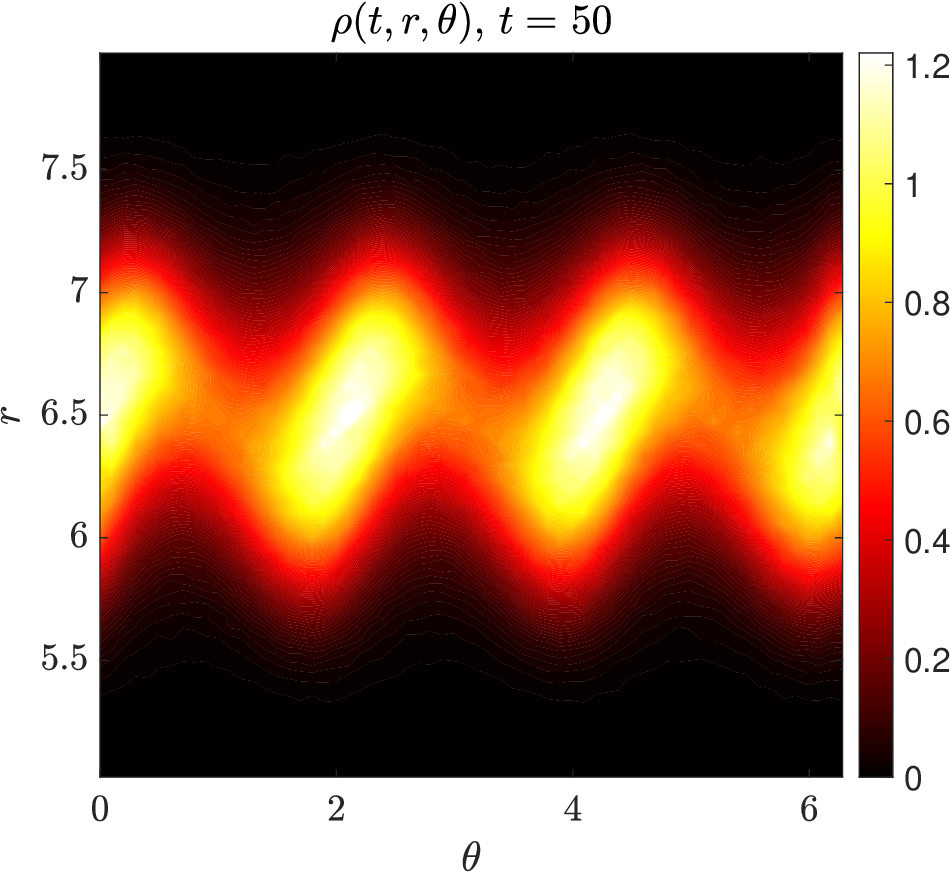}
	\includegraphics[width=0.328\linewidth]{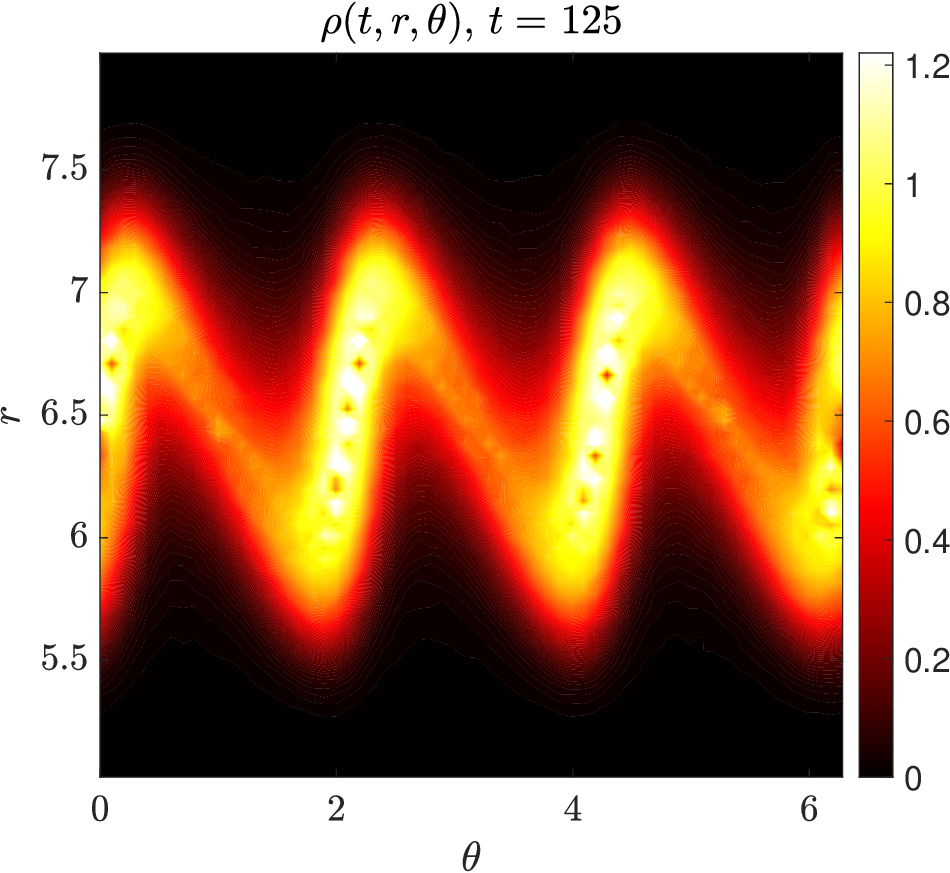}
	\includegraphics[width=0.328\linewidth]{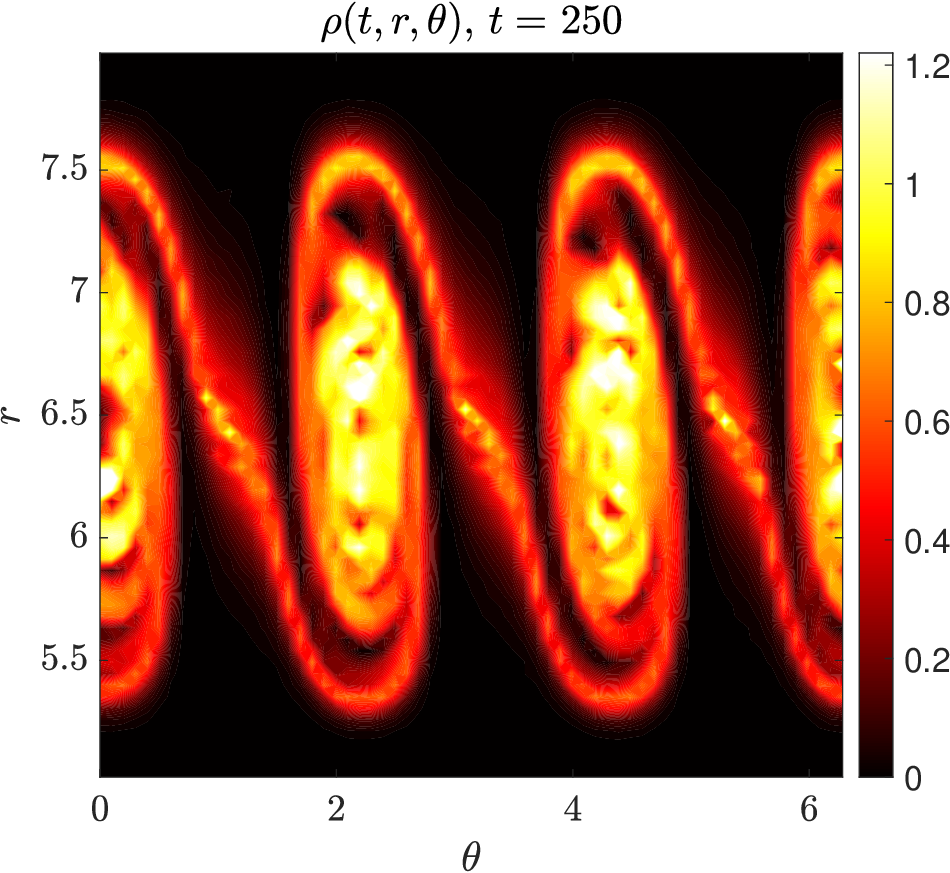}\\
	\includegraphics[width=0.328\linewidth]{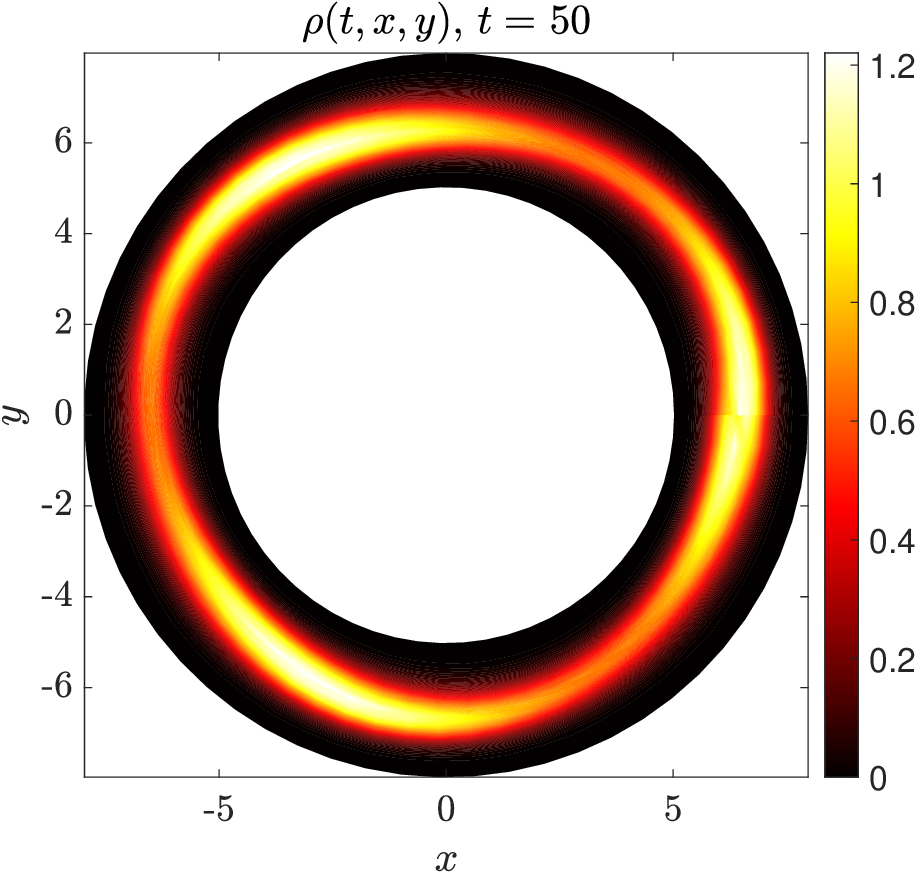}
	\includegraphics[width=0.328\linewidth]{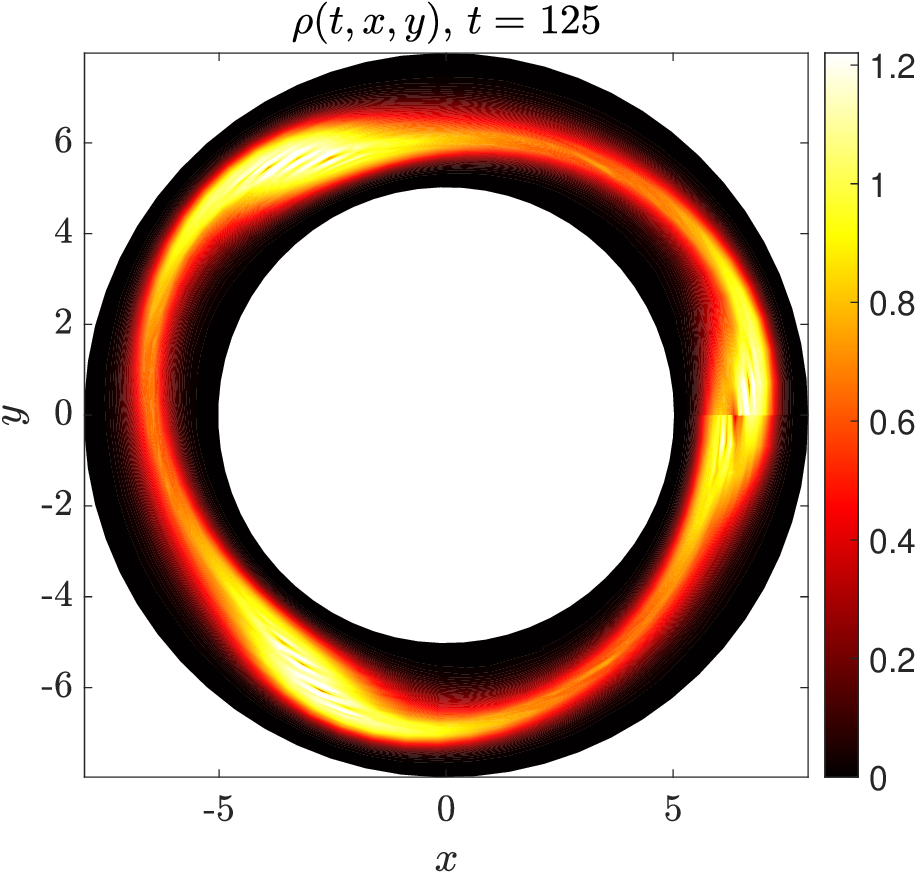}
	\includegraphics[width=0.328\linewidth]{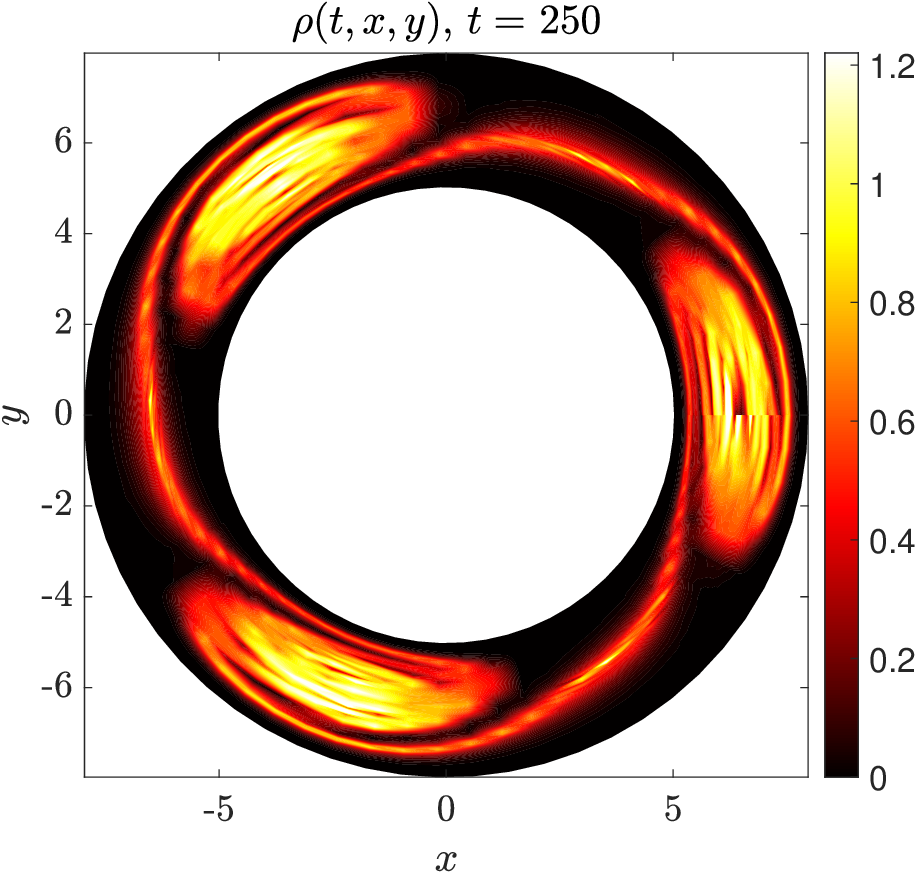}
	\caption{Diocotron instability. Uncontrolled dynamics obtained by setting $B(t,\rr) = 10$. Three snapshots of the dynamics taken at time $t=50$, $t=125$ and $t=200$. First row: polar coordinates. Second row: cartesian coordinates. }
	\label{fig:uncontrolled}
\end{figure}

\begin{figure}[h!]
	\centering
	\includegraphics[width=0.458\linewidth]{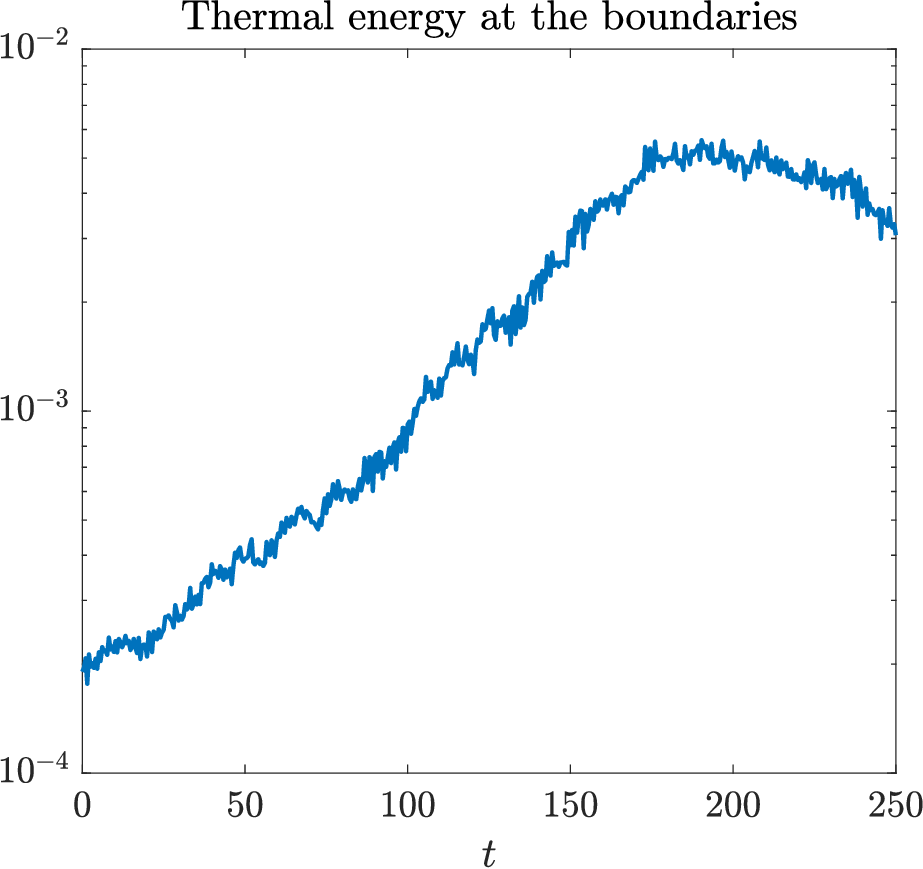}
	\caption{Diocotron instability: uncontrolled dynamics. Thermal energy at the boundaries computed as in \eqref{eq:energy}. }
	\label{fig:uncontrolled_energy}
\end{figure}

We begin by considering the uncontrolled case, setting $B(t,\rr) = 10$.  Figure \ref{fig:uncontrolled} shows three snapshots of the dynamics taken at time $t= 50$, $t=125$ and $t=250$. In the first row, the plot in polar coordinates, in the second row the same plot but in cartesian coordinates. The presence of a strong magnetic field leads to the formation of vortices, as expected. 
In Figure \ref{fig:uncontrolled_energy} the thermal energy at the boundaries computed as in \eqref{eq:energy}. Once that the instability arises, the thermal energy increases. 
\subsection{Test with control}\label{sec:test_control}
\begin{figure}[h!]
	\centering
	\includegraphics[width=0.328\linewidth]{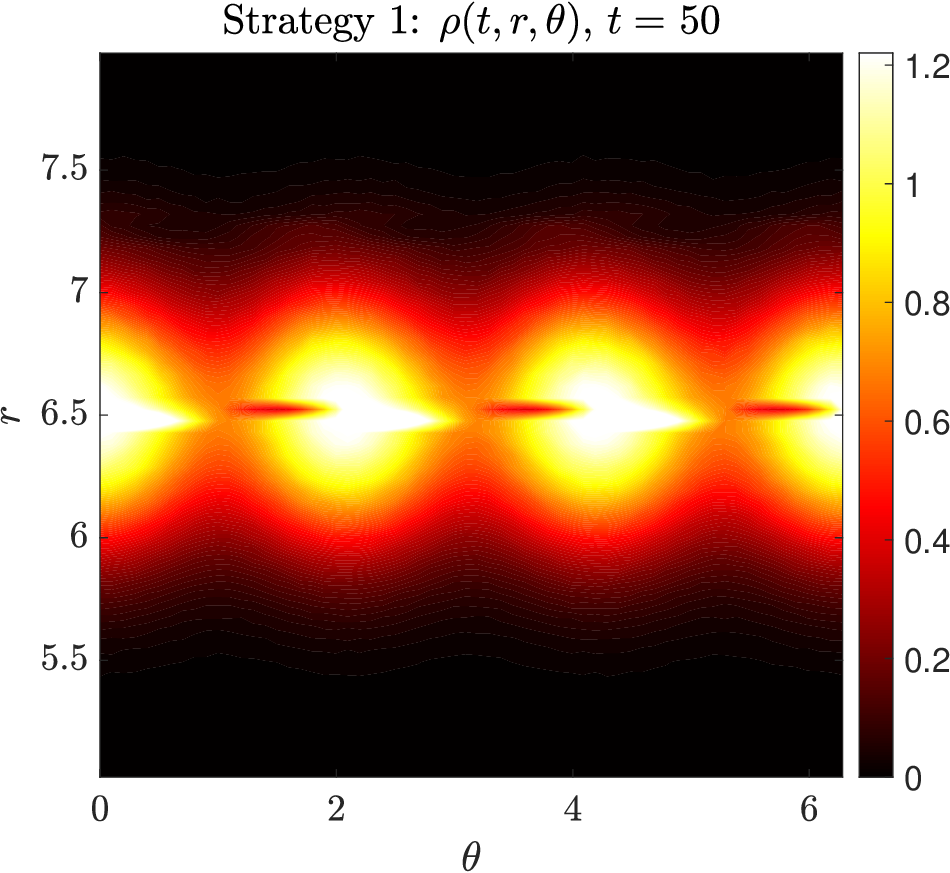}
	\includegraphics[width=0.328\linewidth]{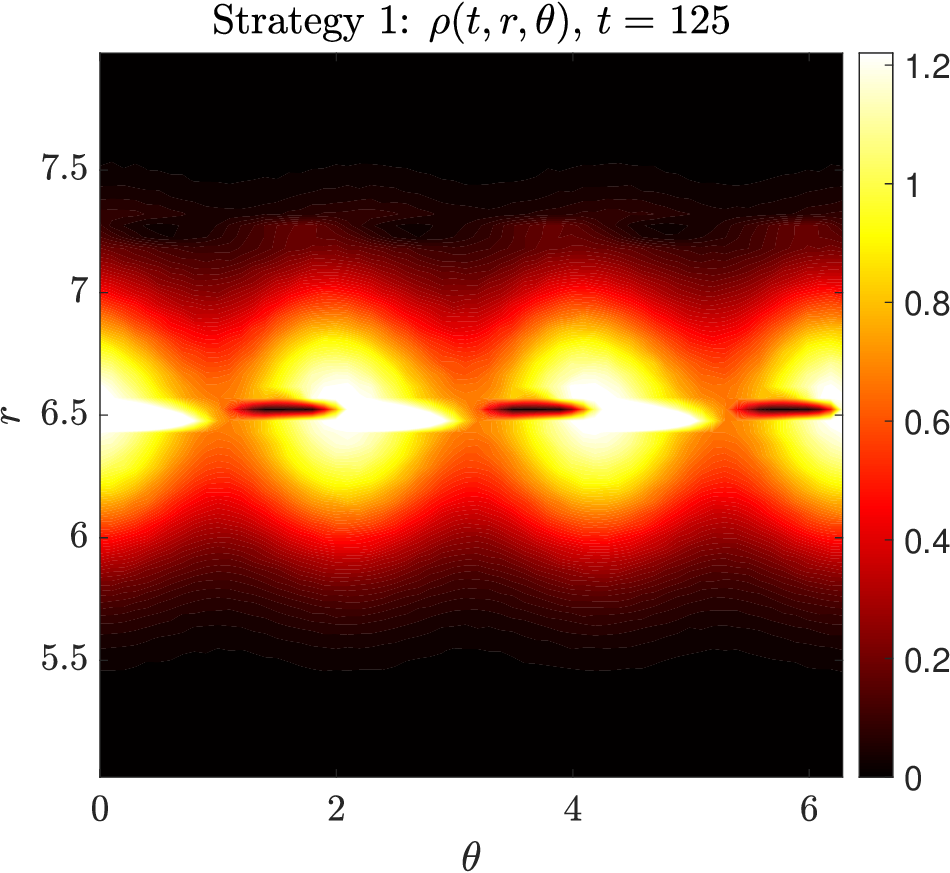}
	\includegraphics[width=0.328\linewidth]{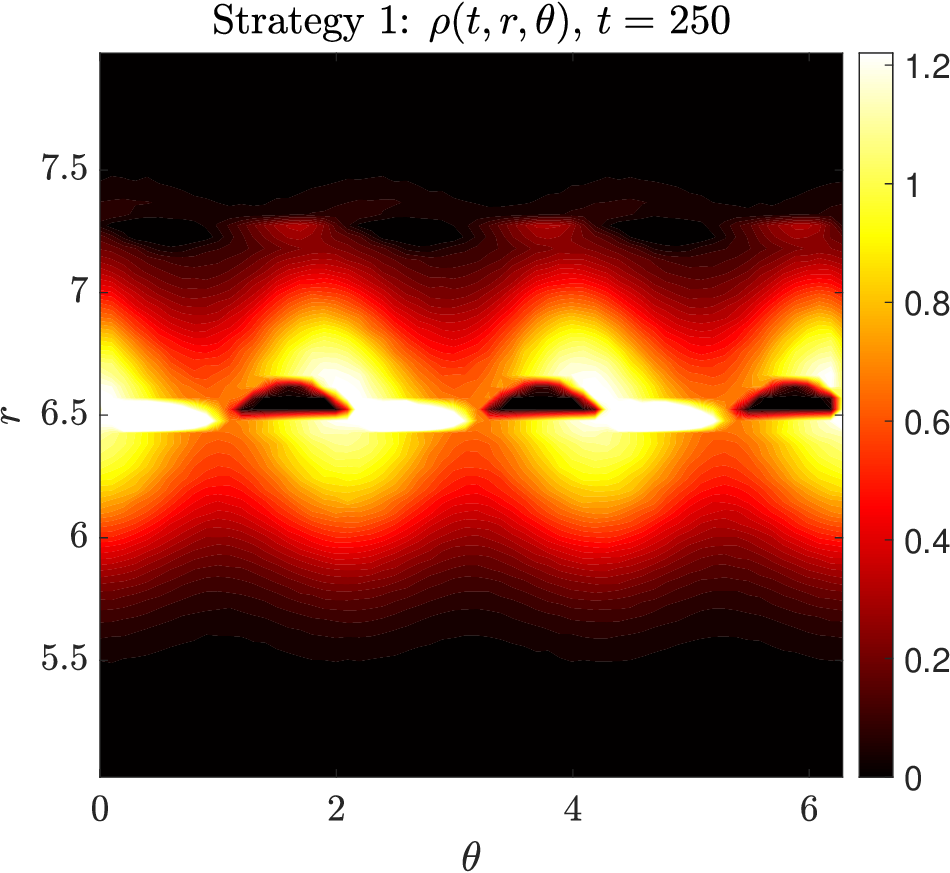}\\
	\includegraphics[width=0.328\linewidth]{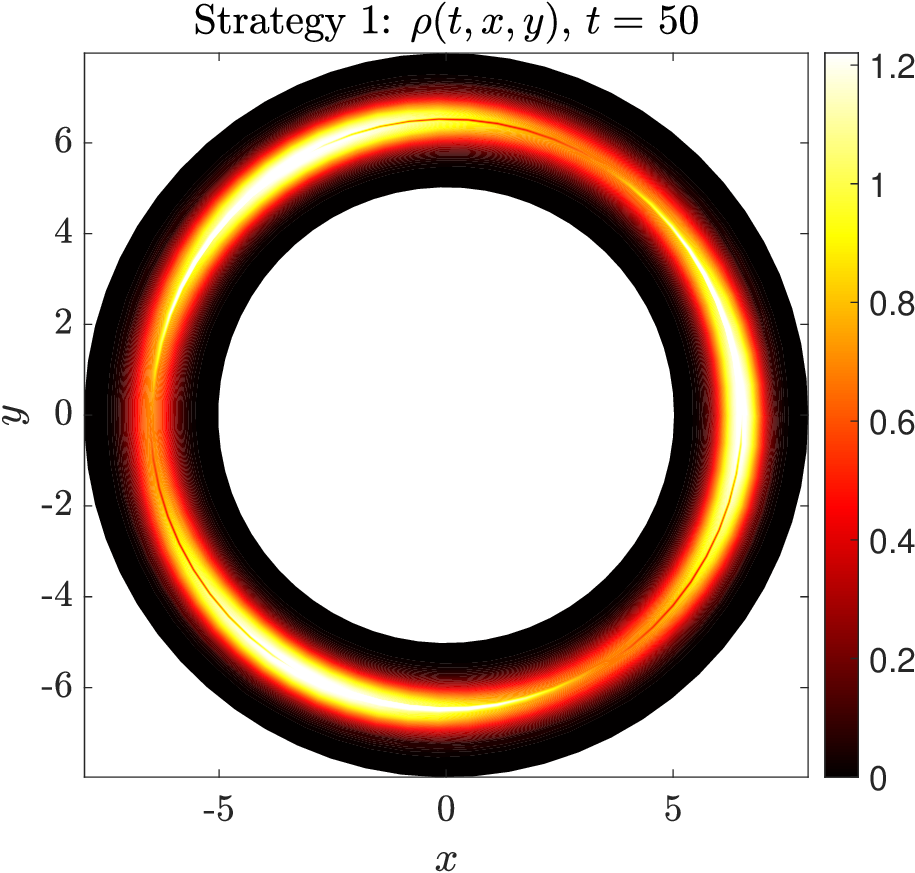}
	\includegraphics[width=0.328\linewidth]{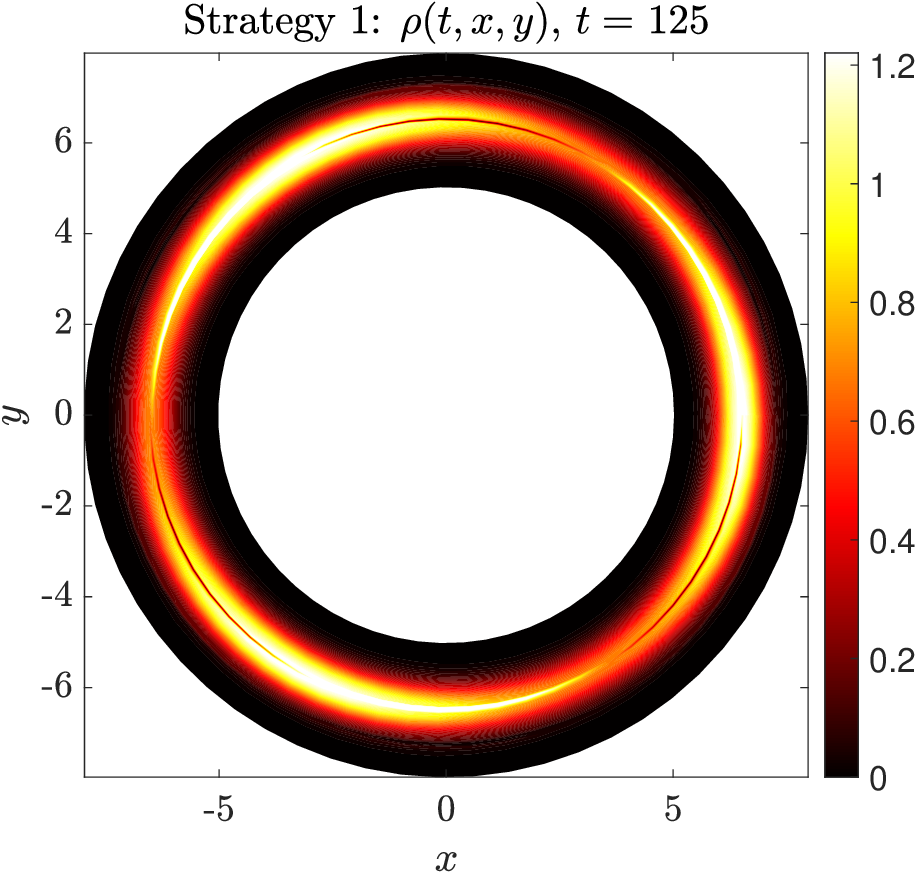}
	\includegraphics[width=0.328\linewidth]{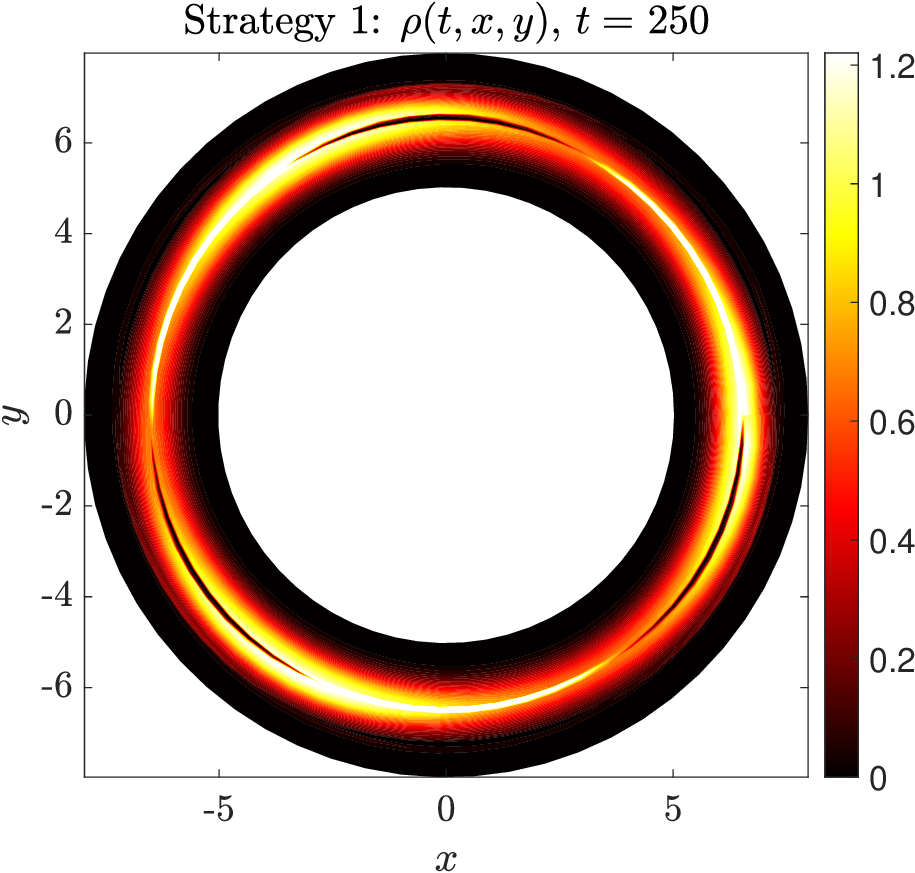}
	\caption{Diocotron instability. Controlled dynamics obtained by setting $B(t,\rr)$ as in \eqref{eq:Bk_strategy1} (strategy one). Three snapshots of the dynamics taken at time $t=50$, $t=125$ and $t=200$. First row: polar coordinates. Second row: cartesian coordinates. }
	\label{fig:strategy1}
\end{figure}

\begin{figure}[h!]
	\centering
	\includegraphics[width=0.328\linewidth]{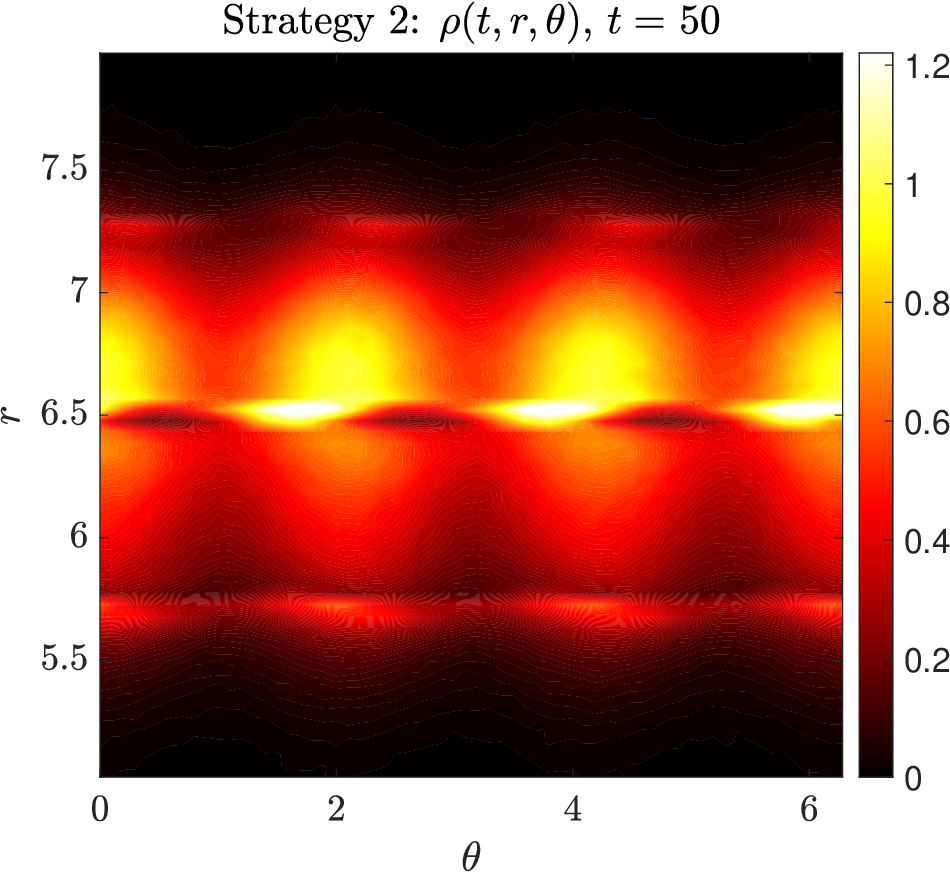}
	\includegraphics[width=0.328\linewidth]{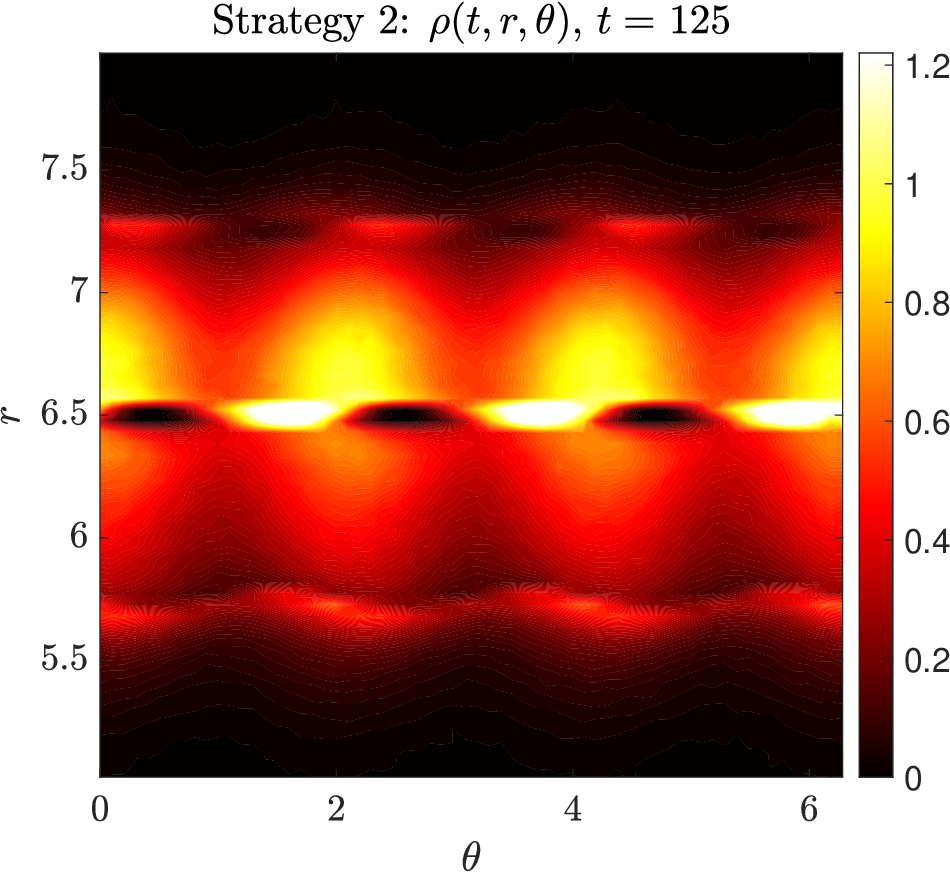}
	\includegraphics[width=0.328\linewidth]{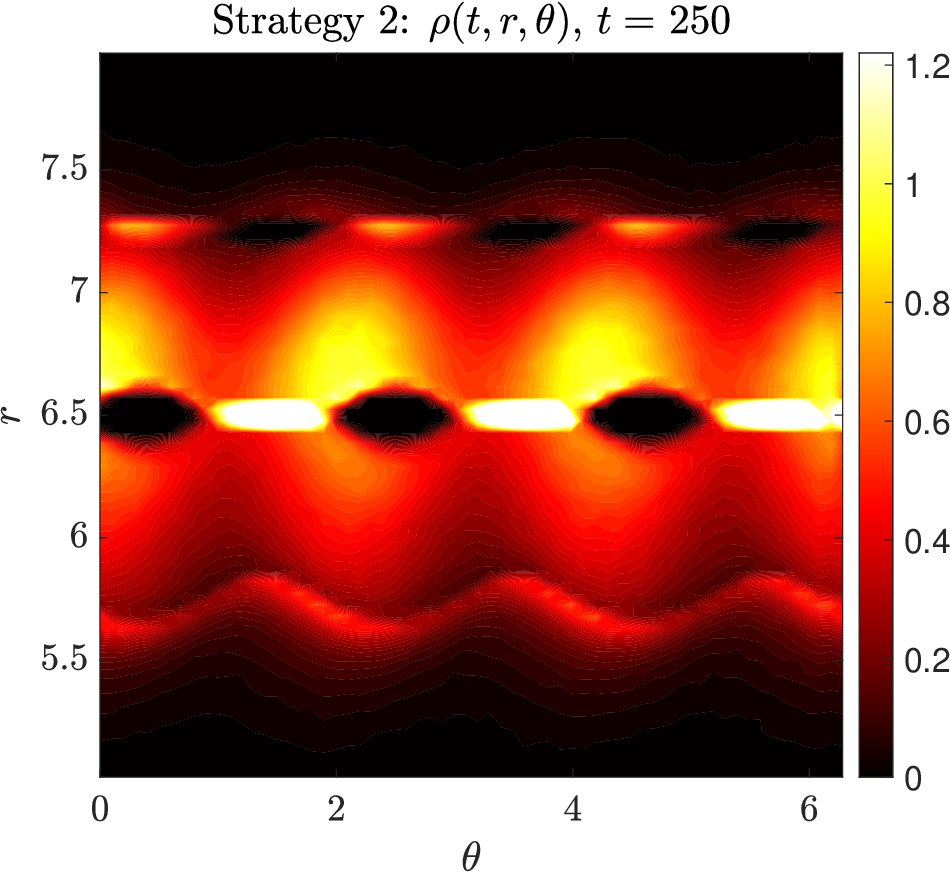}\\
	\includegraphics[width=0.328\linewidth]{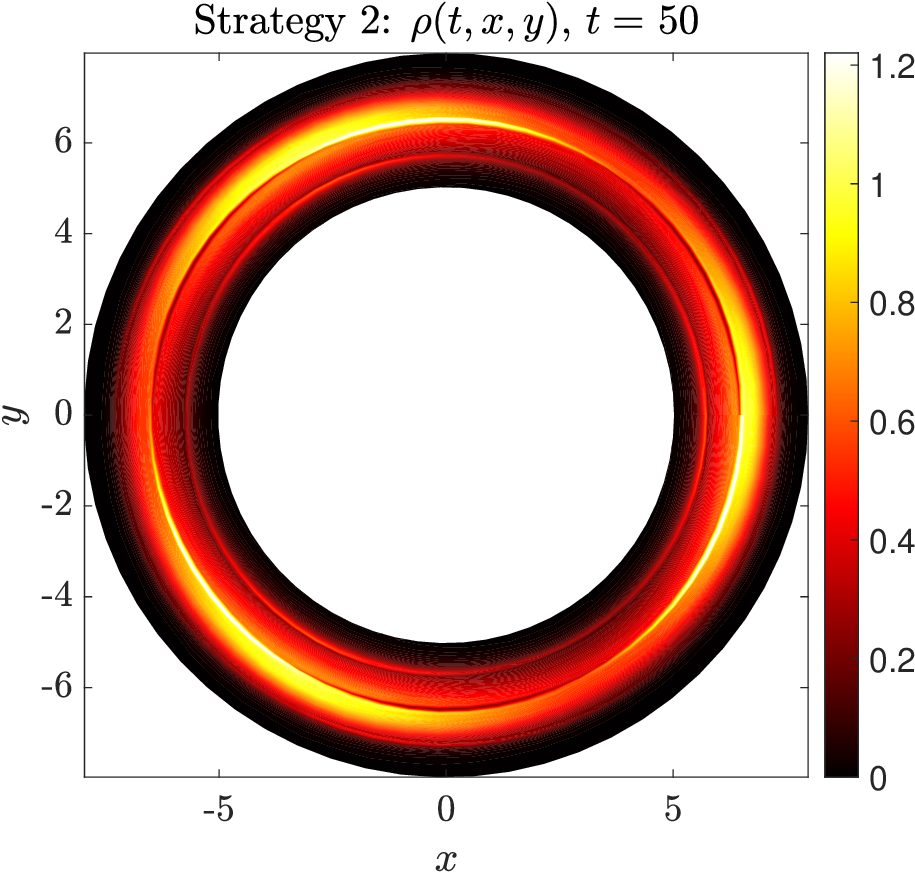}
	\includegraphics[width=0.328\linewidth]{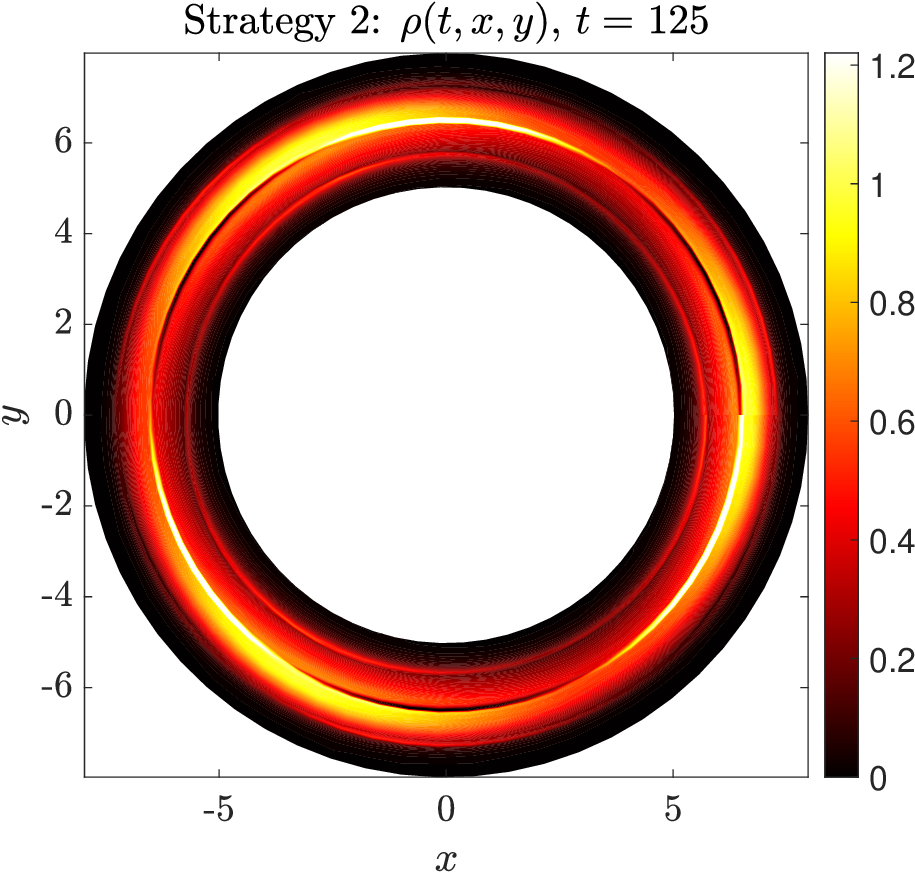}
	\includegraphics[width=0.328\linewidth]{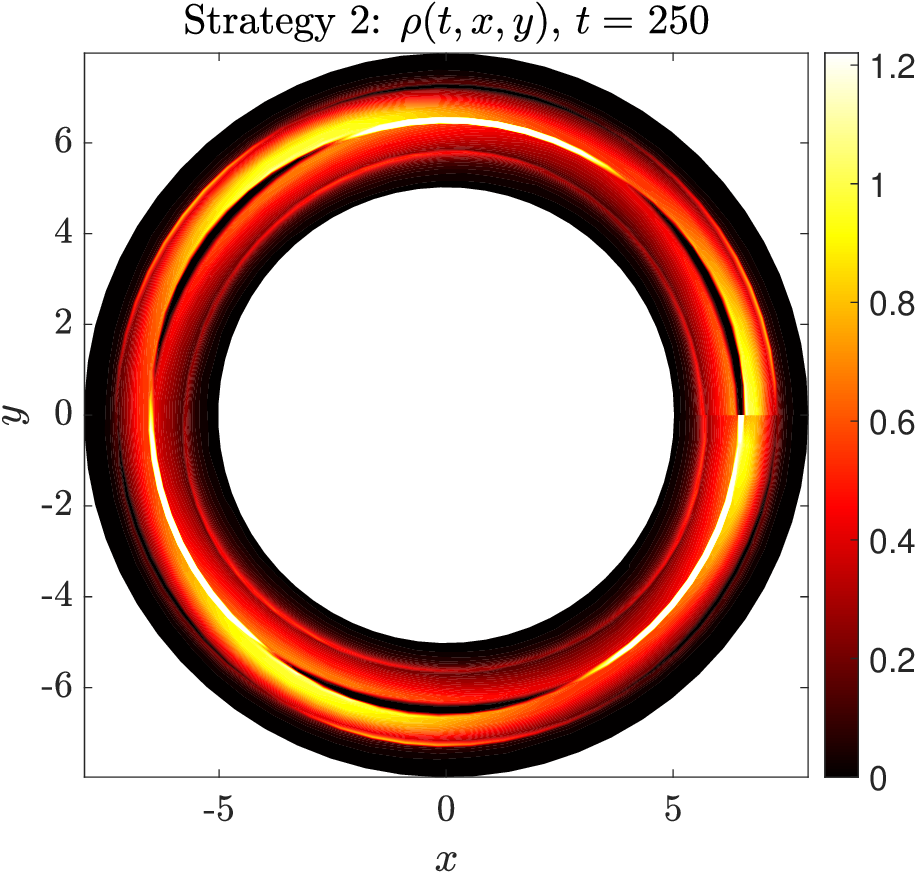}
	\caption{Diocotron instability. Controlled dynamics obtained by setting $B(t,\rr)$ as in \eqref{eq:Bk_strategy2_interp} (strategy two). Three snapshots of the dynamics taken at time $t=50$, $t=125$ and $t=200$. First row: polar coordinates. Second row: cartesian coordinates. }
	\label{fig:strategy2}
\end{figure}

We now test the effectiveness of the control strategies proposed in the previous section. We set $B(t,\rr)$ as in \eqref{eq:Bk_strategy1}-\eqref{eq:Bk_strategy2_interp} with $\alpha_\emph{r} = 100$, $\alpha_\emph{v} = 5$, $\beta_\emph{r} = 10$, $\beta_\emph{v} = 5$, $\gamma=10^{-4}$, $M=100$ and assuming the control takes piecewise constant values on a grid with $N_c=4$ horizontal cell in the $r-\theta$ plane. Figure \ref{fig:strategy1}-\ref{fig:strategy2} show three snapshots of the dynamics at time $t=50$, $t=125$ and $t=250$, referred to the strategy one and two respectively. In the first row we plot the results in polar coordinates, while in the second row in cartesian coordinates.  

\begin{figure}[h!]
	\centering
	\includegraphics[width=0.458\linewidth]{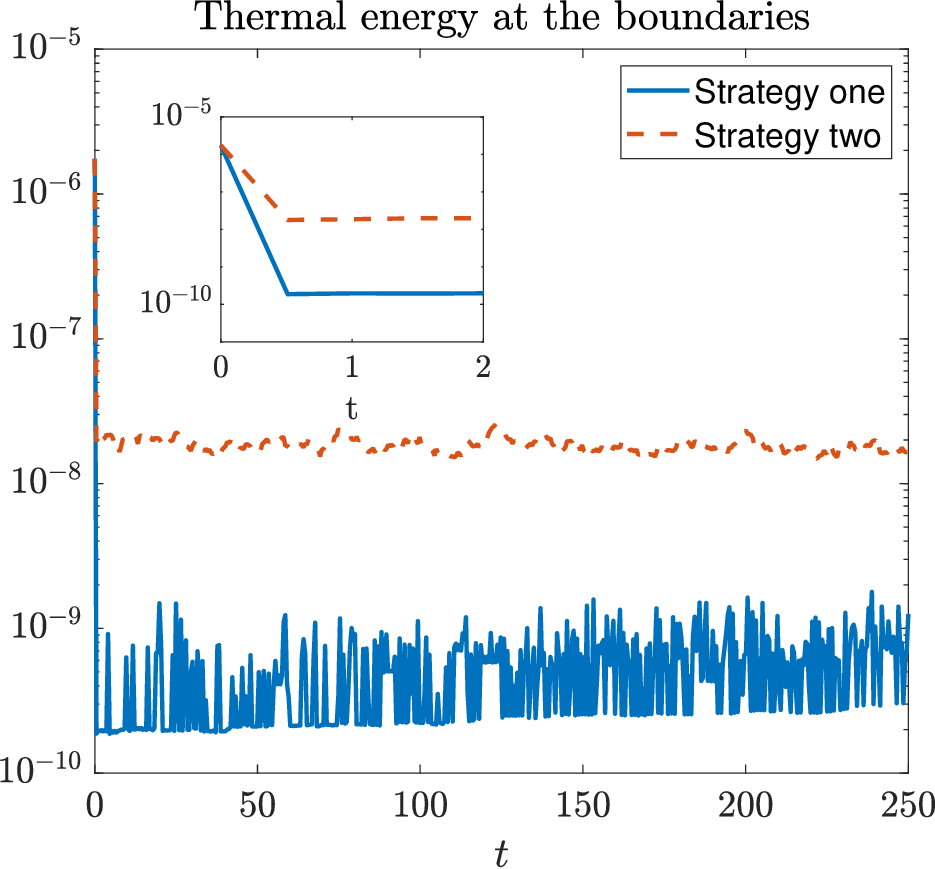} 
	\includegraphics[width=0.45\linewidth]{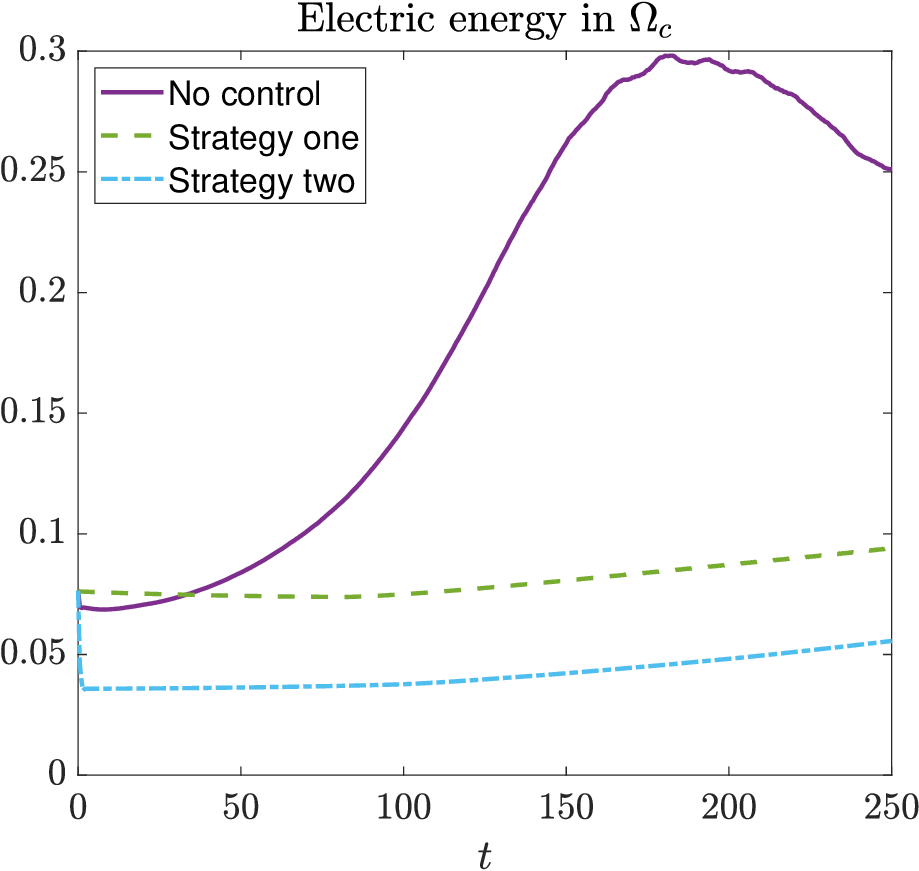} 
	\caption{Diocotron instability. Time evolution of the thermal energy (left) and electric energy (right), computed according to \eqref{eq:energy}–\eqref{eq:electric_energy}, for the case in which $B(t,\rr)$ is defined as in \eqref{eq:Bk_strategy1} (strategy one) and as in \eqref{eq:Bk_strategy2_interp} (strategy two). The inset in the left panel provides a zoomed-in view of the short-time behavior. In the right panel, the electric energy obtained in the absence of control is also shown for comparison. }
	\label{fig:energy}
\end{figure}

\begin{figure}[h!]
	\centering
	\includegraphics[width=0.458\linewidth]{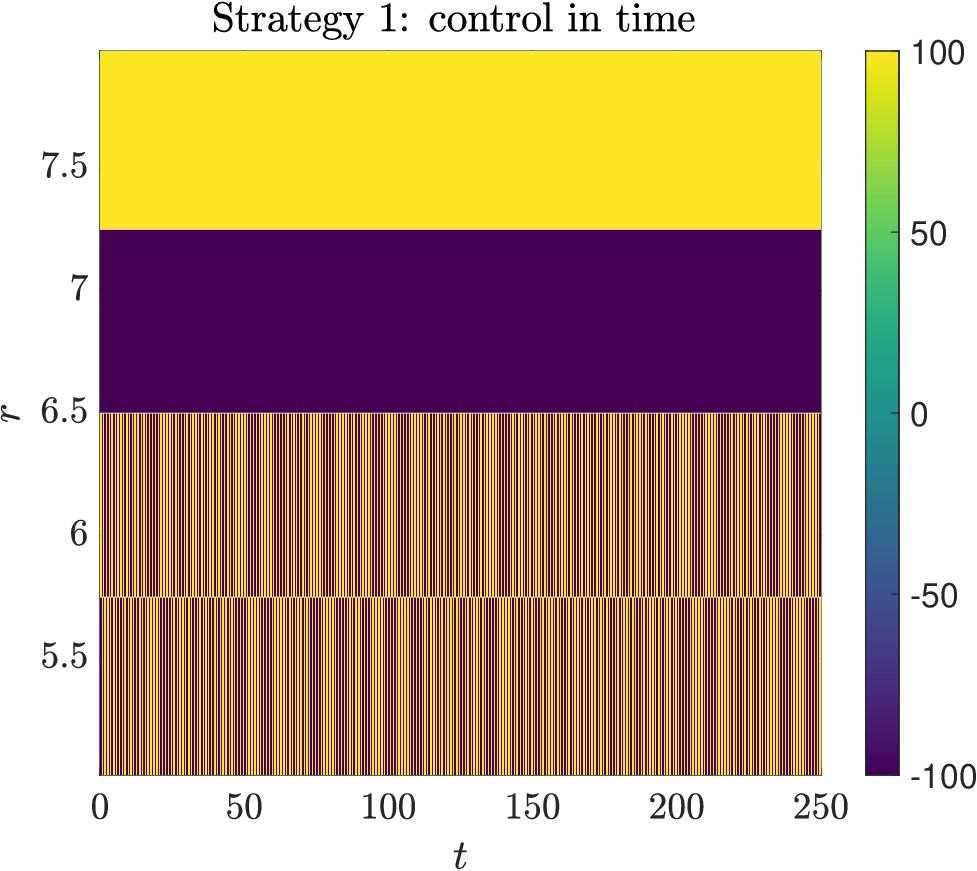}
	\includegraphics[width=0.458\linewidth]{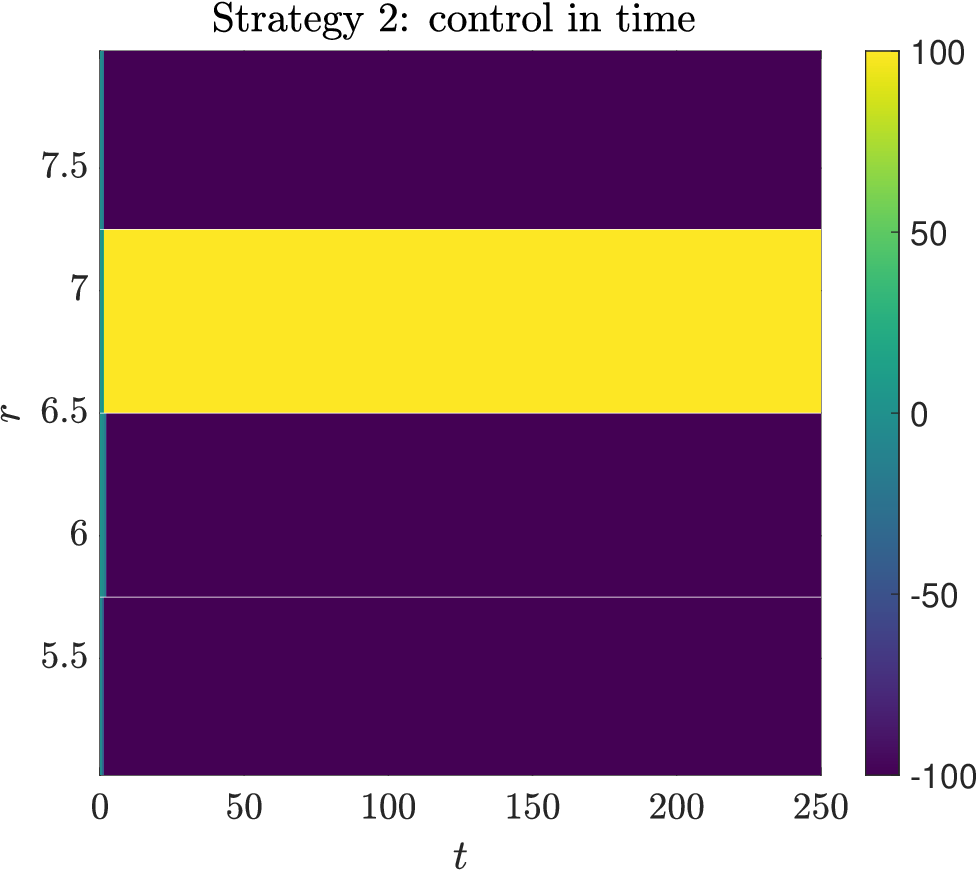}
	\caption{Diocotron instability. Value of the magnetic field in time as in \eqref{eq:Bk_strategy1} (on the left), and as in \eqref{eq:Bk_strategy2_interp} (on the right). }
	\label{fig:control_B}
\end{figure}

Both strategies effectively control the system dynamics, as confirmed by Figure \ref{fig:energy} (left), which depicts the evolution of the thermal energy at the boundaries over time. The embedded plot provides a zoomed-in view of the short-time behavior, highlighting the initial decay of the energy.  Strategy one results in a slightly greater reduction of boundary thermal energy compared to strategy two. The application of the controls induces oscillations in the density, particularly in the center of the domain, where regions of high density alternate with regions of low density. This behavior is primarily due to the piecewise-constant structure of the control, which can confine particles at the interfaces of the fictitious macro-cells $C_k$, as the control exhibits sharp variations between adjacent cells, particularly in the central region of the domain (see Figure \ref{fig:control_B}). Nevertheless, no increase is observed in the energy associated with the self-induced electric field measured at the center of the domain, as defined in \eqref{eq:electric_energy}. We also observe that the electric energy corresponding to strategy two is slightly lower than that obtained with strategy one. This difference may be attributed to the higher particle concentration at the center of the domain produced by strategy one, which leads to a stronger self-induced electric field and, consequently, to an increase in the associated electric energy. This is illustrated in Figure~\ref{fig:energy} (right), where the electric energy in the uncontrolled case is also shown for comparison. 

In Figure \ref{fig:control_B}  the value of the control $B(t,\rr)$ in time, defined as in \eqref{eq:Bk_strategy1} on the left, and as in \eqref{eq:Bk_strategy2_interp} on the right.  
The control obtained with strategy one exhibits significant temporal oscillations in the lower region of the $(r,\theta)$ domain, which serve to guide the particles toward the desired configuration at the center. In contrast, the control from strategy two is more stable over time.   

Table \ref{tab:costs} reports the computational time, in seconds, required for simulations up to the final time $t_f = 250$, both with and without control, comparing the two proposed strategies. In both cases, the inclusion of the control leads to an increase in computational cost. This is due to the additional operations required, namely identifying the position of each particle within the cells to compute their mean velocity and position in strategy one, or interpolating the pointwise values of the control in strategy two.
	Strategy two is slightly more computationally demanding than strategy one, which can be attributed to the higher cost of the interpolation procedure. Overall, strategy one achieves a favourable balance between control effectiveness and computational efficiency.  All the numerical experiments have been run on a Desktop Computer equipped with an Intel(R)
	18 Core(TM) i7-8700 CPU processor and 32GB RAM.
\begin{table}  
	\begin{center}
		\begin{tabular}{|c|c|c|c|}
			\hline
			& No control  & Strategy one & Strategy two \\
			\hline
			Time (seconds)  & 509 & 1112 & 1211 \\
			\hline
		\end{tabular} 
		\caption{Diocotron instability. Computational time in seconds of a simulation up to time $t_f=250$, without control and with control strategy one and strategy two.  }
		\label{tab:costs}
	\end{center} 
\end{table} 
\section{Conclusion} \label{sec:conclusion}
In this work, we have presented an overview on the control strategies introduced in \cite{albi2025instantaneous,albi2025robust} by reformulating the results within the polar coordinates framework. Specifically, we have developed a piecewise control acting on the magnetic field, capable of steering the system toward desired configurations. The first strategy computes the piecewise constant control directly by minimizing a prescribed functional. The second strategy, instead, first determines a pointwise control for each particle and then interpolates these values over the spatial grid to obtain piecewise constant controls. The control problem is formulated over a single time step, resulting in an instantaneous feedback law obtained through a simplified time integration scheme, which is subsequently incorporated into the semi-implicit dynamics.
Both strategies are shown to be effective in reducing thermal energy at the boundaries. Numerical experiments in the two-dimensional polar coordinate framework confirm these results.

Future research will focus on including uncertainty and collisions, providing a more comprehensive framework as in \cite{albi2025robust}. We also plan to tackle the added complexity of the full Maxwell–Vlasov–BGK system and to explore the incorporation of more realistic collision operators, such as the Landau operator, within the proposed control approach. Further directions include developing robust control strategies based on the solution of the optimality conditions and employing neural networks to accurately determine control values. Extensions to more realistic scenarios, including the three-dimensional case in cylindrical coordinates, will also be considered.
\section*{Acknowledges} 
This work has been written within the activities of GNCS group of INdAM
(Italian National Institute of High Mathematics). The research was partially supported by the Italian Ministry of University and Research (MUR) through the PRIN 2020 project (No. 2020JLWP23) "Integrated
Mathematical Approaches to Socio-Epidemiological Dynamics" and funded with the contribution of the Italian Ministry of University and Research (MUR) pursuant to Decree No. 1236 of 01/08/2023 – Fondo Italiano per la Scienza 2022-2023 - Bando FIS 2.
	
\bibliographystyle{abbrv}

\begin{thebibliography}{99.}%
	\bibitem{albi2017mean}
	Giacomo Albi, Young-Pil Choi, Massimo Fornasier, and Dante Kalise, ``Mean field control hierarchy'' \textit{Applied Mathematics \& Optimization}, 76(1), 93-135, 2017.
	\bibitem{albi2025instantaneous}
	Giacomo Albi, Giacomo Dimarco, Federica Ferrarese, and Lorenzo Pareschi,  
	``Instantaneous control strategies for magnetically confined fusion plasma,''  
	\textit{Journal of Computational Physics}, vol.~527, p.~113804, 2025.
	
	\bibitem{albi2025robust}
	Giacomo Albi, Giacomo Dimarco, Federica Ferrarese, and Lorenzo Pareschi,  
	``Robust feedback control of collisional plasma dynamics in presence of uncertainties,'' \textit{Journal of Computational Physics}, vol.~ 546, p.~114512, 2026. 
	
	\bibitem{bartsch2024controlling}
	Jan Bartsch, Patrik Knopf, Stefania Scheurer, and Jörg Weber,  
	Controlling a Vlasov–Poisson plasma by a Particle-in-Cell method based on a Monte Carlo framework,   
	\textit{SIAM Journal on Control and Optimization}, vol.~62, no.~4, pp.~1977--2011, 2024.
	
	\bibitem{belaouar2009asymptotically}
	Radoin Belaouar, Nicolas Crouseilles, Pierre Degond, and Eric Sonnendrücker,  
	An asymptotically stable semi-Lagrangian scheme in the quasi-neutral limit,  
	\textit{Journal of Scientific Computing}, vol.~41, pp.~341--365, 2009.
	
	\bibitem{caprino2012magnetic}
	Stefano Caprino, Giacomo Cavallaro, and Carlo Marchioro,  
	Time evolution of a Vlasov–Poisson plasma with magnetic confinement, 
	\textit{Kinetic and Related Models}, vol.~5, no.~4, pp.~729--742, 2012.
	
	\bibitem{cheng1976integration} 
	Chio-Zong Cheng, and Georg Knorr.  
	The integration of the Vlasov equation in configuration space.  
	\textit{Journal of Computational Physics}, vol.~22, no.~3, p.~330-351, 1976. 
	
	\bibitem{cheng2014discontinuous}
	Yingda Cheng, Irene M. Gamba, Fengyan Li, and Philip J. Morrison,  
	Discontinuous Galerkin methods for the Vlasov–Maxwell equations, 
	\textit{SIAM Journal on Numerical Analysis}, vol.~52, no.~2, pp.~1017--1049, 2014.
	
	\bibitem{chacon2016pic}
	Eugenio Chacón-Golcher, Sorin A. Hirstoaga, and Michael Lutz,  
	Optimization of Particle-in-Cell simulations for Vlasov–Poisson system with strong magnetic field, 
	\textit{ESAIM: Proceedings and Surveys}, vol.~53, pp.~177--190, 2016.
	
	\bibitem{Chac2023}
	Gang Chen and Luis Chacón,  
	An implicit, conservative and asymptotic-preserving electrostatic Particle-In-Cell algorithm for arbitrarily magnetized plasmas in uniform magnetic fields,   
	\textit{Journal of Computational Physics}, vol.~487, Paper No.~112160, 16 pp., 2023.
	
	\bibitem{coughlin2022efficient}
	James Coughlin and Jingwei Hu,  
	Efficient dynamical low-rank approximation for the Vlasov–Ampère–Fokker–Planck system,  
	\textit{Journal of Computational Physics}, vol.~470, p.~111590, 2022.
	
	\bibitem{crouseilles2004highorder}
	Nicolas Crouseilles and Francis Filbet,  
	Numerical approximation of collisional plasmas by high order methods, 
	\textit{Journal of Computational Physics}, vol.~201, no.~2, pp.~546--572, 2004.
	
	\bibitem{crouseilles2004numerical}
	Nicolas Crouseilles and Francis Filbet,  
	Numerical approximation of collisional plasmas by high order methods,  
	\textit{Journal of Computational Physics}, vol.~201, no.~2, pp.~546--572, 2004.
	
	\bibitem{crouseilles2010conservative}
	Nicolas Crouseilles, Michel Mehrenberger, and Eric Sonnendrücker,  
	Conservative semi-Lagrangian schemes for Vlasov equations,   
	\textit{Journal of Computational Physics}, vol.~229, no.~6, pp.~1927--1953, 2010.
	
	\bibitem{crouseilles2016multiscale}
	Nicolas Crouseilles, Giacomo Dimarco, and Marie-Hélène Vignal,  
	Multiscale schemes for the BGK–Vlasov–Poisson system in the quasi-neutral and fluid limits. Stability analysis and first order schemes,  
	\textit{Multiscale Modeling \& Simulation}, vol.~14, no.~1, pp.~65--95, 2016.
	\bibitem{crouseilles2013semi}
	Nicolas Crouseilles, P. Glanc, S. Hirstoaga, E. Madaule, M. Mehrenberger, and J. Pétri,  
	Semi-Lagrangian simulations on polar grids: from diocotron instability to ITG turbulence test cases,  
	in \textit{4th International Workshop on the Theory and Applications of the Vlasov Equation (VLASOVIA 2013)}, 2013.
	
	\bibitem{degond2013ap}
	Pierre Degond,  
	Asymptotic-preserving schemes for fluid models of plasmas,   
	in \textit{Numerical Models for Fusion}, Panorama \& Synthèses, vol.~39/40, pp.~1--90, Soc. Math. France, Paris, 2013.
	
	\bibitem{degond2017ap}
	Pierre Degond and Franck Deluzet,  
	Asymptotic-preserving methods and multiscale models for plasma physics, 
	\textit{Journal of Computational Physics}, vol.~336, pp.~429--457, 2017.
	
	\bibitem{dimarco2010dsmc}
	Giacomo Dimarco, Ronald E. Caflisch, and Lorenzo Pareschi,  
	Direct simulation Monte Carlo schemes for Coulomb interactions in plasmas,  
	\textit{Communications in Applied and Industrial Mathematics}, vol.~1, no.~1, pp.~1--72, 2010.
	
	\bibitem{dimarco2015numerical}
	Giacomo Dimarco, Qin Li, Lorenzo Pareschi, and Bokai Yan,  
	Numerical methods for plasma physics in collisional regimes, 
	\textit{Journal of Plasma Physics}, vol.~81, no.~1, 2015.
	
	\bibitem{einkemmer2024vlasov}
	Lukas Einkemmer, Qi Li, Lin Wang, and Yunan Yang,  
	Suppressing instability in a Vlasov–Poisson system by an external electric field through constrained optimization, 
	\textit{Journal of Computational Physics}, vol.~498, p.~112662, 2024.
	
	\bibitem{einkemmer2025control}
	Lukas Einkemmer, Qi Li, Clément Mouhot, and Yue Yunan,  
	Control of instability in a Vlasov–Poisson system through an external electric field,  
	\textit{Journal of Computational Physics}, vol.~530, p.~113904, 2025.
	
	\bibitem{fasoli2016computational}
	Ambrogio Fasoli, Stephan Brunner, Wilfred A. Cooper, John P. Graves, Paolo Ricci, Olivier Sauter, and Laurent Villard,  
	Computational challenges in magnetic-confinement fusion physics, 
	\textit{Nature Physics}, vol.~12, no.~5, pp.~411--423, 2016.
	
	\bibitem{filbet2003numerical} 
	Francis Filbet, and Eric Sonnendrücker.  
	Numerical methods for the Vlasov equation.  
	\textit{Numerical Mathematics and Advanced Applications: Proceedings of ENUMATH 2001 the 4th European Conference on Numerical Mathematics and Advanced Applications Ischia, July 2001}. Milano: Springer Milan, 2003.
	
	\bibitem{filbet2016pic}
	Francis Filbet and Luís M. Rodrigues,  
	Asymptotically stable Particle-in-Cell methods for the Vlasov–Poisson system with a strong external magnetic field, 
	\textit{SIAM Journal on Numerical Analysis}, vol.~54, no.~2, pp.~1120--1146, 2016.
	
	\bibitem{filbet2017ap}
	Francis Filbet and Luís M. Rodrigues,  
	Asymptotically preserving Particle-in-Cell methods for inhomogeneous strongly magnetized plasmas,  
	\textit{SIAM Journal on Numerical Analysis}, vol.~55, no.~5, pp.~2416--2443, 2017.
	
	\bibitem{filbet2018numerical} 
	Francis Filbet, and Chang Yang.  
	Numerical Simulations to the Vlasov-Poisson System with a Strong Magnetic Field.  
	\textit{HAL,} 2018.
	
	\bibitem{garabedian2003computational}
	Paul R. Garabedian,  
	Computational mathematics and physics of fusion reactors,
	\textit{Proceedings of the National Academy of Sciences}, vol.~100, no.~24, pp.~13741--13745, 2003.
	
	\bibitem{ghizzo1993eulerian}
	Alain Ghizzo, Pierre Bertrand, Magdi Shoucri, Eric Fijalkow, and Marc R. Feix,  
	An Eulerian code for the study of the drift-kinetic Vlasov equation,   
	\textit{Journal of Computational Physics}, vol.~108, no.~1, pp.~105--121, 1993.
	
	\bibitem{grandgirard2013gyrokinetic}
	Virginie Grandgirard and Yves Sarazin,  
	Gyrokinetic simulations of magnetic fusion plasmas,   
	\textit{Panoramas et Synthèses}, no.~39–40, pp.~91--176, 2013.
	
	\bibitem{gu2022hamiltonian}
	Aiguo Gu, Yu He, and Yu Sun,  
	Hamiltonian Particle-in-Cell methods for Vlasov–Poisson equations, 
	\textit{Journal of Computational Physics}, vol.~467, p.~111472, 2022.
	
	\bibitem{hairer2018boris}
	Earl Hairer and Christian Lubich,  
	Energy behaviour of the Boris method for charged-particle dynamics,  
	\textit{BIT Numerical Mathematics}, vol.~58, pp.~969--979, 2018.
	
	\bibitem{han-kwan2010tokamak}
	Daniel Han-Kwan,  
	On the confinement of a Tokamak plasma, 
	\textit{SIAM Journal on Mathematical Analysis}, vol.~42, no.~6, pp.~2337--2367, 2010.
	
	\bibitem{imbert2024introduction}
	Lise-Marie Imbert-Gérard, Elizabeth J. Paul, and Adelle M. Wright,  
	An Introduction to Stellarators: From Magnetic Fields to Symmetries and Optimization,  
	\textit{SIAM}, 2024.
	
	\bibitem{knopf2020optimal}
	Patrik Knopf and Jörg Weber,  
	Optimal control of a Vlasov–Poisson plasma by fixed magnetic field coils, 
	\textit{Applied Mathematics \& Optimization}, vol.~81, pp.~961--988, 2020.
	
	\bibitem{deluzet2023parallelization}
	Franck Deluzet, Guillaume Fubiani, Loïc Garrigues, Clément Guillet, and Jérôme Narski,  
	Efficient parallelization for 3D-3V sparse grid Particle-in-Cell: Shared memory architectures,  
	\textit{Journal of Computational Physics}, vol.~480, p.~112022, 2023.
	
	\bibitem{russo2009semilagrangian}
	Giovanni Russo and Francis Filbet,  
	Semilagrangian schemes applied to moving boundary problems for the BGK model of rarefied gas dynamics,   
	\textit{Kinetic and Related Models}, vol.~2, no.~1, pp.~231--250, 2009.
	
	\bibitem{sonnendrucker1999semi}
	Eric Sonnendrücker, Jean Roche, Pierre Bertrand, and Alain Ghizzo,  
	The semi-Lagrangian method for the numerical resolution of the Vlasov equation, 
	\textit{Journal of Computational Physics}, vol.~149, no.~2, pp.~201--220, 1999.
	
	\bibitem{sonnendrucker1999semilagrangian}
	Eric Sonnendrücker, Jean Roche, Pierre Bertrand, and Alain Ghizzo,  
	The semi-Lagrangian method for the numerical resolution of the Vlasov equation, 
	\textit{Journal of Computational Physics}, vol.~149, no.~2, pp.~201--220, 1999.
	
	\bibitem{spitzer1958stellarator}
	Lyman Spitzer,  
	The stellarator concept,   
	\textit{Physics of Fluids}, vol.~1, no.~4, p.~253, 1958.
	
	\bibitem{taitano2013implicit}
	William T. Taitano, David A. Knoll, Luis Chacón, and Gang Chen,  
	Development of a consistent and stable fully implicit moment method for Vlasov–Ampère Particle-In-Cell (PIC) system, 
	\textit{SIAM Journal on Scientific Computing}, vol.~35, no.~5, pp.~S126--S149, 2013.
	
	\bibitem{weber2021optimal}
	Jörg Weber,  
	Optimal control of the two-dimensional Vlasov–Maxwell system,   
	\textit{ESAIM: Control, Optimisation and Calculus of Variations}, vol.~27, p.~S19, 2021.
	
	\bibitem{yang2014hermite}
	Cheng Yang and Francis Filbet,  
	Conservative and non-conservative methods based on Hermite weighted essentially non-oscillatory reconstruction for Vlasov equations, 
	\textit{Journal of Computational Physics}, vol.~279, pp.~18--36, 2014.
	\bibitem{valentini2005numerical}
	Francesco Valentini, Pierluigi Veltri, and André Mangeney,  
	A numerical scheme for the integration of the Vlasov--Poisson system of equations, in the magnetized case,  
	\textit{Journal of Computational Physics}, vol.~210, no.~2, pp.~730--751, 2005.
	
\end{thebibliography}

\end{document}